\documentclass[journal,12pt,onecolumn,draftclsnofoot,]{IEEEtran}
\usepackage{amsmath,amsfonts}
\usepackage{algorithmic}
\usepackage{algorithm}
\usepackage{array}
\usepackage[caption=false,font=normalsize,labelfont=sf,textfont=sf]{subfig}
\usepackage{textcomp}
\usepackage{stfloats}
\usepackage{url}
\usepackage{verbatim}
\usepackage{graphicx}
\usepackage{cite}
\begin{document}
	
	\title{Efficient Service Differentiation and Energy Management in Hybrid WiFi/LiFi Networks}
	
	\author{Asim~Ihsan, Muhammad~Asif, Hossein Safi, Iman Tavakkolnia, \IEEEmembership{Member,~IEEE}, and Harald Haas,~\IEEEmembership{Fellow,~IEEE}
		\thanks{Asim Ihsan, Hossein Safi, Iman Tavakkolnia, and Harald Haas  are with LiFi Research and Development Centre (LRDC) , Electrical Engineering Division, Department of Engineering, University of Cambridge, Cambridge, UK, (Emails: {ai422, hs905, it360, huh21}@cam.ac.uk).}
		\thanks{Muhammad Asif is with the School of Computer Science and
			Communication Engineering, Jiangsu University, Zhenjiang, China. (Email:
			masif@ujs.edu.cn).}
		\thanks{This work was presented in part as an extended abstract at the 52nd Wireless World Research Forum (WWRF52), held at King's College London, UK, from September 10–12, 2024.}
		\thanks{This work was supported in part by the Project REASON, a UK Government funded Project
			through the Future Open Networks Research Challenge (FONRC) sponsored by the Department of
			Science Innovation and Technology (DSIT), and in part by the Engineering and Physical Sciences
			Research Council (EPSRC) through the Platform for Driving Ultimate Connectivity (TITAN),
			EP/X04047X/1 and EP/Y037243/1.}
		
	}
	% The paper headers
	%\markboth{Journal of \LaTeX\ Class Files,~Vol.~14, No.~8, August~2021}%
	%{Shell \MakeLowercase{\textit{et al.}}: A Sample Article Using IEEEtran.cls for IEEE Journals}
	
	%\IEEEpubid{0000--0000/00\$00.00~\copyright~2021 IEEE}
	% Remember, if you use this you must call \IEEEpubidadjcol in the second
	% column for its text to clear the IEEEpubid mark.
	
	\maketitle
	
	\begin{abstract} 
		 In this paper, we propose an innovative predict-and-optimize algorithm designed for hybrid WiFi/LiFi networks, aiming to achieve service differentiation while maximizing energy efficiency (EE). The proposed framework utilizes multi-access technology real-time intelligent controller (mATRIC) to dynamically predict the appropriate network slice for each user based on historically monitored key performance indicators (KPIs). This prediction is facilitated by a deep learning model trained using the resilient backpropagation algorithm, with training conducted on KPIs data at the universal non-real time  RAN intelligent controller (non-RT RIC). This trained model enables real-time slice selection by mATRIC. In the subsequent phase, the algorithm focuses on optimizing EE of hybrid network as a function of precoding vectors for the predicted slices by employing techniques from sequential convex approximation and the inner approximation method. We introduce novel approximations to convert non-convex objective functions and constraints into convex forms, and develop an iterative algorithm to achieve sub-optimal solutions. Additionally, the EE maximization problem, ensures alignment with end-to-end latency requirements. It also addresses the various constraints inherent to hybrid systems, such as input signal limitations for LiFi LEDs, data rate restrictions, and power budget considerations. Simulation results validate the effectiveness of the proposed algorithm, demonstrating significant improvements in EE while ensuring service differentiation within hybrid network environments.
	\end{abstract}
	
	
	\section{Introduction and motivation}
	The exponential growth in mobile data traffic and the diverse requirements of emerging applications necessitates innovative solutions in wireless communication systems. The advent of 5G technology and subsequent generations has highlighted the potential of network slicing. Network slicing allows for the differentiation of services by segmenting the network into distinct slices, each adapted to specific application needs \cite{5G Network Slicing}. This capability supports the heterogeneous needs of enhanced mobile broadband (eMBB), ultra-reliable low-latency communication (URLLC), and massive machine-type communications (mMTC), while striving to balance performance with energy efficiency. Despite its advantages in managing diverse service requirements, network slicing poses challenges, particularly related to energy consumption and operational costs, raising concerns about sustainability \cite{EE Slicing}. As network operators aim to deliver a broad range of services, integrating energy-efficient strategies becomes essential.
	
	Open radio access network (O-RAN) is an initiative to develop an open, disaggregated RAN architecture that promotes interoperability and flexibility by standardizing interfaces and components. This approach allows for greater innovation and cost-efficiency, enabling seamless integration of different radio technologies into a cohesive network ecosystem. O-RAN will unify various radio technologies, including non-3GPP technologies, Wi-Fi and LiFi, to satisfy the diverse requirements of future wireless networks and create a more integrated and versatile communication infrastructure \cite{5G-ClARITY}. For control and optimization of RAN functions across multi-access technologies, including all wired and wireless technologies, the realizing enabling architectures and solutions for open networks (REASON) project plays a crucial role. REASON is dedicated to advancing the O-RAN initiative by developing multi-access technology real-time intelligent controller (mATRIC), an innovative framework designed to manage and optimize a diverse array of access technologies \cite{mATRIC}. It acts as a central component, receiving policy directives from the universal non-real time  RAN intelligent controller (non-RT RIC) and executing real-time control and optimization of multi-access technologies. By utilizing collected data, mATRIC ensures dynamic network performance and adaptability and supports various xApps. Here, xApps are microservice based application that are executed on the mATRIC to enhance network functionality and efficiency in real-time.
	
	In O-RAN, non-3GPP access networks, like WiFi and LiFi, can be integrated using non-3GPP interworking functions. Both WiFi and LiFi networks can connect through the same interworking function, but each network uses its own separate IP subnetwork. To fully capitalize on the benefits of mATRIC for hybrid WiFi/LiFi network, it is crucial to extend network slicing to accommodate these hybrid technologies and applications within the RAN ecosystem \cite{Hybrid Slice}. Consequently, supporting network slicing for each technology is essential. For Wi-Fi and Li-Fi networks, slicing can be effectively applied to address their unique requirements, allowing for the creation of tailored virtual slices that optimize performance and resource allocation. This approach ensures that each network can meet specific needs, such as high-density Wi-Fi environments or high-speed Li-Fi communications, enhancing overall efficiency and service quality \cite{WiFi Slice, WiFi Slice 2, LiFi Slice}. By applying network slicing to a hybrid setup, it is possible to optimize the use of both WiFi and LiFi resources in a coordinated manner through mATRIC. mATRIC can dynamically allocate resources among different access technologies and network slices, ensuring optimal performance for eMBB, URLLC, and mMTC applications.
	
	To effectively support the diverse service requirements of eMBB, URLLC, mMTC in a hybrid WiFi/LiFi network, it is crucial to consider the specific transmission needs of each service type. Each of these services has distinct performance and reliability requirements that impact the design of the network and the optimization strategies employed. For URLLC, which demands extremely low latency and high reliability, the transmission duration plays a critical role. Given that the end-to-end latency requirement for URLLC is approximately 1 millisecond \cite{URLLC end-to-end 1,URLLC end-to-end 2}, a transmission duration of 0.05 milliseconds \cite{URLLC time 1,URLLC time 2}, is sufficient to meet this stringent requirement \cite{URLLC time 3}. Finite block length transmission is crucial here due to the short block lengths involved, which affect error rates. Traditional Shannon capacity, applicable to longer block lengths, does not accurately reflect achievable rates for URLLC. Instead, the Polyanskiy bound \cite{FBL Information Theory}, considering finite block lengths, provides a more accurate measure of achievable rates, which are lower than those predicted by Shannon capacity. This necessitates algorithms that handle finite block lengths effectively to ensure high reliability and meet URLLC's stringent performance requirements. For mMTC, which manages numerous devices with intermittent communication needs, finite block length transmission is particularly advantageous. While mMTC's latency requirements are less stringent than those of URLLC, short block lengths are still valuable due to the nature of mMTC’s communication. This technology often involves transmitting short payloads, where ensuring reliability and efficiency is crucial. Finite block length transmission enables the system to handle these brief communication intervals effectively, maintaining reliable performance, thus aligning well with the demands of mMTC \cite{mMTC FBL}. Dealing with the finite block length constraints for URLLC and mMTC, alongside traditional large block length transmission for eMBB, presents a complex challenge for resource allocation. This scenario introduces novel optimization problems that require advanced strategies for efficient resource management \cite{O-RAN Network slicing}.
	
    In hybrid WiFi/LiFi networks, most studies have contributed towards key issues such as EE \cite{EE Hybrid 1,EE Hybrid 2}, load balancing \cite{Load Balancing 1,Load Balancing 2}, handover \cite{Handover 1}, user association \cite{User association 1, User association 2}, and access point selection \cite{Access Point selection 1}. However, a notable gap exists in integrating service differentiation with EE, an aspect that remains unexplored in hybrid WiFi/LiFi networks. Addressing this gap, our work focuses on developing strategies that not only balance performance and energy consumption but also meet diverse service requirements of different application, thereby fully exploiting the potential of hybrid networks. Our major contribution is summarized as follows
	\begin{itemize}
		\item A novel predict-and-optimize algorithm is introduced for hybrid WiFi/LiFi networks to achieve service differentiation while maximizing EE. The framework first employs mATRIC to dynamically predict the network slice for each user based on monitored KPIs, using a deep learning model trained in MATLAB with the resilient backpropagation algorithm. The model can be trained on monitored KPIs at the universal non-RT RIC which then can be used for slice selection for users in real-time by mATRIC. In the subsequent stage, EE of the hybrid network is optimized for the selected slices of users by applying techniques from sequential convex approximation and the inner approximation method. We develop new approximations to convexify non-convex objective functions and constraints, and an iterative algorithm is used to find sub-optimal solutions. The EE maximization process also incorporates end-to-end latency requirements, and is addressing constraints for hybrid network such as input signal limitations for LiFi LEDs, data rate restrictions, and power budget considerations.  Performance is evaluated through simulations, demonstrating improved EE and effective service differentiation in hybrid network environments.
		\item The proposed algorithm’s worst-case computational complexity is computed using Big O notation and is observed to be low, ensuring both practical efficiency and scalability. This makes the algorithm suitable for real-world deployment in hybrid WiFi/LiFi networks, as both the neural network prediction and EE optimization phases handle complex tasks with manageable computational demands.
		\item Numerical results are presented, demonstrating the fast convergence of the proposed algorithm and its effectiveness in optimizing transmit precoding vectors for hybrid WiFi/LiFi networks with RAN slicing support. This results in enhanced EE while maintaining effective service differentiation.
		\item In our study, we compared the performance of our proposed optimization algorithm for EE maximization as function of precoding vectors in hybrid networks with conventional precoding methods, such as Zero-Forcing (ZF) and Maximum Ratio Transmission (MRT). The results indicate that our algorithm significantly outperforms these traditional methods in terms of EE while effectively ensuring service differentiation. Additionally, we assessed the performance of a fully hybrid LiFi/WiFi system versus a WiFi-only system. The hybrid system, which makes use of both WiFi and LiFi transmitters, demonstrates superior EE compared to the WiFi-only system. This highlights the benefits of integrating LiFi with WiFi, not only in optimizing resource use but also in enhancing service differentiation and overall network performance.
	\end{itemize}
The paper is further proceeds as follows: The system model and problem formulation are discussed in Section II. Section III covers the solution of formulated problem, and Section IV provides the numerical results along with its related discussions. At last, Section V delivers concluding remarks and suggestions for future research. The definitions of the different symbols used in our research are provided in Table I.  

\begin{table}[h]
	\centering
	\caption{The table of different symbols and their definitions}
	\begin{tabular}{|c|l|} % Added vertical lines
		\hline % Top horizontal line
		\textbf{Symbol} & \textbf{Definition} \\ 
		\hline % Line between header and body
		$K$  & Total number of users in the system \\ 
		\hline
		$M$  & The number of antennas at the WiFi transmitter \\ 
		\hline
		$L$  & The number of LEDs at the LiFi transmitter \\
		\hline
		$s$  & Type of service, where $s \in \{eMBB, URLLC, mMTC\}$  \\ 
		\hline
		$\mathbf{h}_{k,s}^\text{WiFi}$  & The channel link between the WiFi transmitter and the $k^{th}$ user in slice $s$   \\ 
		\hline
		$\mathbf{h}_{k,s}^\text{LiFi}$ & The channel link between the LiFi transmitter and the $k^{th}$ user in slice $s$  \\ 
		\hline
		$\mathbf{f}_{k,s}^\text{WiFi}$  & WiFi precoder for the $k^{th}$ user in slice $s$  \\ 
		\hline
		$\mathbf{f}_{k,s}^\text{LiFi}$  & LiFi precoder for the $k^{th}$ user in slice $s$ \\ 
		\hline
		$(\sigma_{k,s}^\text{WiFi})^2$  & Variance of AWGN for the $k^{th}$ user in WiFi slice $s$  \\ 
		\hline
		$(\sigma_{k,s}^\text{LiFi})^2$  & Variance of AWGN for the $k^{th}$ user in LiFi slice $s$  \\ 
		\hline
		$\gamma_{k,s}^\text{WiFi}$  & SINR of the $k^{th}$ user in WiFi slice $s$  \\ 
		\hline
		$\gamma_{k,s}^\text{LiFi}$  & SINR of the $k^{th}$ user in LiFi slice $s$ \\ 
		\hline
		$R_{k,s}^\text{WiFi}$  & Achievable rate of the $k^{th}$ user in WiFi slice $s$  \\ 
		\hline
		$R_{k,s}^\text{LiFi}$  & Achievable rate of the $k^{th}$ user in LiFi slice $s$ \\ 
		\hline
		${\text{S}_{k,s}^\text{WiFi}}$ & Shanon term of the $k^{th}$ user in WiFi slice $s$ \\
		\hline
		${\text{V}_{k,s}^\text{WiFi}}$ & Channel dispersion term of the $k^{th}$ user in WiFi slice $s$  \\
		\hline
		${\text{S}_{k,s}^\text{LiFi}}$ & Shanon term of the $k^{th}$ user in LiFi slice $s$ \\
		\hline
		${\text{V}_{k,s}^\text{LiFi}}$ & Channel dispersion term of the $k^{th}$ user in LiFi slice $s$  \\
		\hline
		$P_{\text{WiFi}}$ & Power consumption of the WiFi transmit signal\\
		\hline
		$P_{\text{LiFi}}$ & Power consumption of the LiFi transmit signal\\
		\hline 
		$T_{k,s}$ &  The total latency for the transmission to the $k^{th}$ user in slice $s$\\
		\hline
	\end{tabular}
	\label{tab:symbols}
\end{table}

	\section{SYSTEM MODEL AND PROBLEM FORMULATION}
	
	In this section, we explore a hybrid WiFi/LiFi system, as illustrated in Fig. 1, which comprises multi LED LiFi transmitter and a multi–antenna WiFi transmitter that serves $K$  users, with each user equipped with both a photodiode and an antenna. Each user is indexed by $k$ such that $k \in \{1,2,\cdots,K\}$. Moreover, the system model considers the size of the room, which can be expressed as $L \times W \times H$ cubic meters, where $L$ is the length, $W$ is the width, and $H$ is the height. For this model, the room dimensions are 10 meters in length, 10 meters in width, and 5 meters in height. The WiFi transmitter is placed at the center of the room at coordinates$(W/2,L/2)$. The ceiling of the room is divided into a grid pattern, with $I$ intervals along both the $X$ and $Y$ axes, resulting in a total of $I \times I$ LED positions. Each LED is positioned at a height of $H-1$ meters above the floor, ensuring a uniform distribution throughout the room. The specific positions of the LEDs are determined by dividing the room width and length into equal segments and placing the LEDs at these intervals. Additionally, the $K$ users are randomly distributed throughout the room, with their locations modeled using a binomial point process (BPP).
	
	\begin{figure*}
		\centering
		\includegraphics[width=160mm]{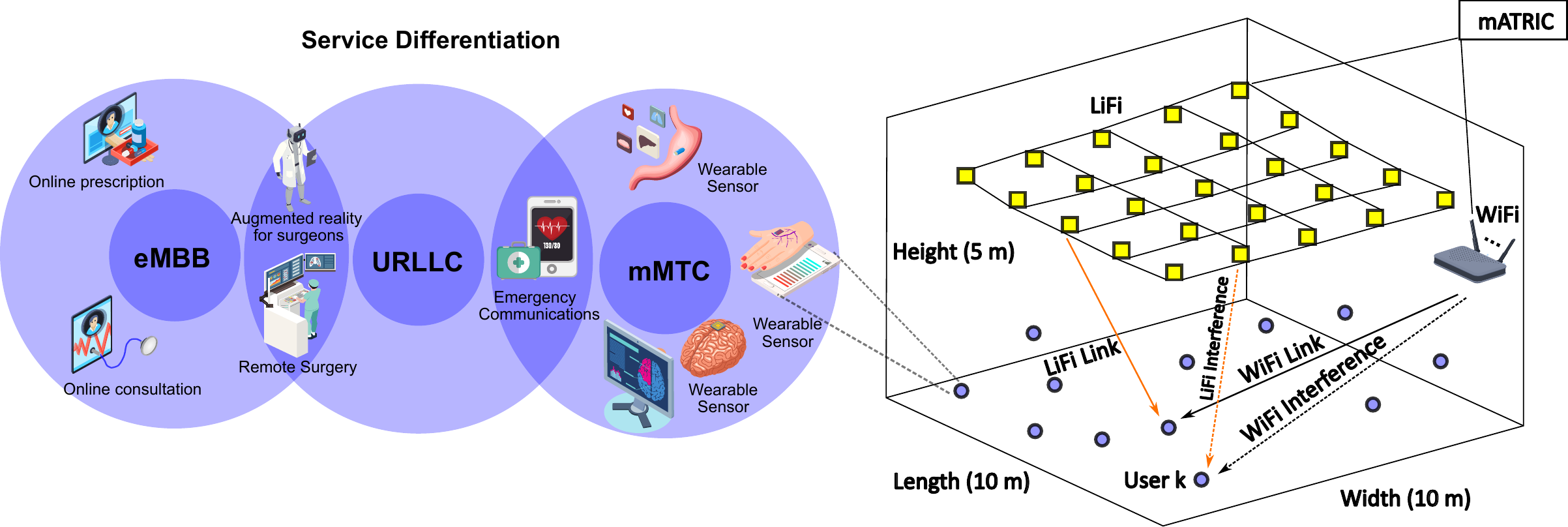}
		\caption{Illustration of System Model}
		\label{FIG:2}
	\end{figure*}
	
	In our analysis, we consider three types of slices: $s \in \{eMBB, URLLC, mMTC\}$ for Hybrid WiFi/LiFi network \cite{LiFi Slice,WiFi Slice}. These slices enable the service differentiation customized to specific needs. mATRIC is tasked with the efficient allocation of network slices and the maximization of EE within a hybrid network environment. Slice selection is driven by KPIs collected by mATRIC from various network nodes such as latency, reliability, supported technologies, use case category, and temporal information from the network infrastructure. Machine learning models is trained on historical KPI data by univeral non-real time (Non-RT) intelligent controller, which can be employed in real time through mATRIC to automate slice selection, ensuring that the network meets the diverse demands of applications and users effectively.

	For WiFi channel, the model for the channel from the WiFi transmitter to the $k^\text{th}$ user at slice $s$ incorporates both channel fading and path loss. It is given by \cite{Access Point selection 1}
	
	\begin{equation}
		\boldsymbol{h}_{k,s}=\Bigg[\sqrt{\frac{\delta}{1+\delta}}e^{j\frac{\pi}{4}}+\sqrt{\frac{1}{1+\delta}}\boldsymbol{g}_{k,s}\Bigg]10^{\frac{-PL[\text{dB}]}{20}},
		\label{eq:1}
	\end{equation}where, $\boldsymbol{g}_{k,s} \in \mathbb{C}^M$, and $M$ is the number of antennas at the WiFi transmitter. The Ricean factor $\delta$ is is set to 1 for distances up to the breakpoint distance $d_{BP}$, and it drops to 0 for distances beyond $d_{BP}$. Additionally, the path loss $PL[\text{dB}]$  is described as
\begin{equation}
	PL[\text{dB}]= \begin{cases} 
		L(d_{k,s})+\Omega_{in}, & d_{k,s}  \leq  d_{BP}, \\
		L(d_{k,s})+35\log_{10}\Big[\frac{d_{k,s}}{d_{BP}}\Big]+\Omega_{out}, & d_{k,s}  >  d_{BP},
	\end{cases},
	\label{eq:2}
\end{equation}
where the free space path loss for WiFi is expressed as
	\begin{equation}
	L(d_{k,s})=20\log_{10}(d_{k,s})+20\log_{10}(f_c)-147.5,
	\label{eq:3}
\end{equation}
with $f_c$ as a central carrier frequency. The shadow fading is modeled as Gaussian random variables, where $\Omega_{in}\sim \mathcal{N}(0, \sigma_{\text{in}}^2)$ for distances within $d_{BP}$, and $\Omega_{out}\sim \mathcal{N}(0, \sigma_{\text{out}}^2)$ for distances beyond $d_{BP}$. For the LiFi channel in our analysis, we extend our modeling approach to include the transmission from a multi-LED LiFi transmitter to the $k^{th}$ user at slice $s$. For LiFi channels, only Line-of-Sight (LoS) links are taken into account, as they carry the majority of the transmitted power, while the power contribution from non-Line-of-Sight (NLoS) links is relatively minimal \cite{LiFi Channel}. The locations of the $l^{th}$ LED and the $k^{th}$ user in a three-dimensional Cartesian coordinate system are represented by $(x_l^\text{LED}, y_l^\text{LED}, z_l^\text{LED})$ and $(x_k^\text{user}, y_k^\text{user}, z_k^\text{user})$ respectively. The distance between the $l^{th}$ LED and the $k^{th}$ user, considering only the LoS path, is given by
	\begin{equation}
		d_{l,k} = \sqrt{(x_l^\text{LED} - x_k^\text{user})^2 + (y_l^\text{LED} - y_k^\text{user})^2 + (z_l^\text{LED} - z_k^\text{user})^2}.
		\label{eq:4}
	\end{equation}
	It is assumed that the transceivers' normal vectors are aligned perpendicularly to the ground plane. Consequently, the angles of irradiance $\Theta_{l,k}$ and incidence $\phi_{l,k}$ for each Line-of-Sight (LoS) component of the channel gain are defined as follows \cite{Surface Congiguration}
	\begin{equation}
		\Theta_{l,k} = \phi_{l,k} = \arccos \left( \frac{z_l^\text{LED} - z_k^\text{user}}{d_{l,k}} \right).
		\label{eq:5}
	\end{equation}
	Now, the LoS channel gain between $l^{th}$ LED and $k^{th}$ user at slice $s$ can be modeled through Lambertian radiant formulation as described below \cite{LiFi Channel Gain}
	\begin{equation}
		h_{l,k,s} =
		\begin{cases} 
			\frac{(m+1) A_D G_f R_D}{2 \pi d_{l,k}^2} \cos^m (\Theta_{l,k}) \cos (\phi_{l,k}) g (\phi_{l,k}), & \text{if } 0 \leq \phi_{l,k} \leq \phi_\text{FoV}, \\
			0, & \text{otherwise}.
		\end{cases}
		\label{eq:6}
	\end{equation}
	The Lambertian index $m$ is determined by the transmitter's half intensity radiation angle $\psi_{1/2}$ as $m = \frac{-1}{\log_2 \left( \cos (\psi_{1/2}) \right)}$. Variables $A_D$, $G_f$, and $R_D$ represent the physical area of the detector, the gain of the optical filter, and the responsivity of the detector, respectively. $\phi_\text{FoV}$ signifies the field of view (FoV) of the detector. The optical concentrator gain $g(\phi_{l,k})$ is defined as
	\begin{equation}
		g(\phi_{l,k}) =
		\begin{cases} 
			\frac{r^2}{\sin^2 (\phi_\text{FoV})}, & \text{if } 0 \leq \phi_{l,k} \leq \phi_\text{FoV}, \\
			0, & \text{otherwise},
		\end{cases}
		\label{eq:7}
	\end{equation}
	where $r$ represents the refractive index, typically ranging between 1 and 2 for visible light.
	
	In LiFi, the detection of signals is based on counting the received photons, which necessitates that the visible light signals be real-valued and non-negative. The transmitted signals for all users are represented by the vector $\mathbf{a}_{s}^\text{LiFi} = \left[ a_{1,s}^\text{LiFi}, a_{2,s}^\text{LiFi}, \ldots, a_{K,s}^\text{LiFi} \right]^T \in \mathbb{R}^{K \times 1}$, and $\mathbf{F}_s^\text{LiFi} = \left[ \mathbf{f}_{1,s}^\text{LiFi}, \mathbf{f}_{2,s}^\text{LiFi}, \cdots, \mathbf{f}_{K,s}^\text{LiFi} \right] \in \mathbb{R}^{L \times K}$ as the corresponding precoding matrix, with $\mathbf{f}_{k,s}^\text{LiFi}$ is the precoding vector for $a_{k,s}^{LiFi}$, which is defined as $\mathbf{f}_{k,s}^\text{LiFi} =  \left[ f_{1,k,s}^\text{LiFi}, f_{2,k,s}^\text{LiFi}, \cdots, f_{L,k,s}^\text{LiFi} \right]^\top \in \mathbb{R}_+^{L \times 1}$. Without loss of generality, these signals are normalized such that $\mathbb{E}\{a_{k,s}^\text{LiFi}\} = 0$ and $\mathbb{E}\{(a_{k,s}^\text{LiFi})^2\} = 1$. Hence, the real-valued input electric signal at the $l^{th}$ LED at the slice $s$ is obtained by summing all the weighted symbols for the $K$ users, given by
	\begin{equation}
		x_{l,s}^\text{LiFi} = \sum_{k=1}^K f_{l,k,s}^\text{LiFi} \, a_{k,s}^\text{LiFi}.
		\label{eq:8}
	\end{equation}
	However,  $x_{l,s}^\text{LiFi}$  cannot be used directly as the input signal to the $l^\text{th}$ LED at slice $s$ because it can sometimes be negative. To ensure proper illumination, a DC bias  $x_{\text{DC},l}^s$  must be added to $x_{l,s}^\text{LiFi}$ to generate a non-negative drive current $I_{l,s}$ for the LED. Additionally, this drive current must be limited to a maximum threshold $I_{max,s}$ to guarantee the normal operation of the LEDs. This can be expressed as
	\begin{equation}
		I_{l,s} = x_{l,s}^\text{LiFi} + x_{\text{DC},l}^s \geq 0\leq I_{\text{max},s}, \quad \forall l.
		\label{eq:9}
	\end{equation}
		To ensure Eq. (9), the following constraint should be satisfied \cite{AC power LiFi}
	\begin{equation}
		\sum_{k=1}^K \left| f_{l,k,s}^\text{LiFi} \right| = \min \left( x_{\text{DC},l}^s, I_{\text{max},s} - x_{\text{DC},l}^s \right).
		\label{eq:10}
	\end{equation} 
	%\begin{equation}
		%I_{l,s} \leq I_{\text{max},s}.
	%\end{equation}
	 $\mathbf{f}_{k,s}^\text{LiFi} = \left[ f_{1,k,s}^\text{LiFi}, f_{2,k,s}^\text{LiFi}, \cdots, f_{L,k,s}^\text{LiFi} \right]^\top \in \mathbb{R}_+^{L \times 1}$ represent the real-valued precoding vector for the $k^{th}$ user from $L$ multi LED LiFi transmitter, and $\mathbf{h}_{k,s}^\text{LiFi}=[h_{1,k,s}^\text{LiFi}, h_{2,k,s}^\text{LiFi}, \ldots, h_{L,k,s}^\text{LiFi}]^\top \in \mathbb{R}_+^{L \times 1}$ denote the channel gain vector for the $k^{th}$ user from multi LED LiFi transmitter. After eliminating the DC component, the received signal at the $k^{th}$ user at the slice $s$ can be expressed as
	
	%Let  $a_{k,s}^\text{LiFi}$ represent the real-valued symbol intended for transmission from a multi-LED LiFi transmitter to the $k^\text{th}$ user at slice $s$, which is assumed to be selected from an MPAM constellation and normalized within the range $[-1, 1]$ \cite{LiFi Precoding}. Then, $\mathbf{a}_s^\text{LiFi}=\left[ a_{1,s}^\text{LiFi}, a_{2,s}^\text{LiFi}, \cdots, a_{K,s}^\text{LiFi} \right] \in \mathbb{R}^{K \times 1}$ represents the symbol vector of all $K$ users at the slice $s$,

	%From Eq. (6) and (7)
	%\begin{equation}
	%	0 \leq I_{l,s} \leq I_{\text{max},s}.
	%\end{equation}
	%Given that  $a_{k,s}^\text{LiFi} \in [-1, 1]$, we determine that the dynamic range of $x_{l,s}$ is
	%\begin{equation}
	%	-\sum_{k=1}^K \left| f_{l,k,s}^\text{LiFi} \right| + x_{\text{DC},l}^s \leq x_{l,s} \leq \sum_{k=1}^K \left| f_{l,k,s}^\text{LiFi} \right| + x_{\text{DC},l}^s.
	%\end{equation}

	\begin{equation}
		y_{k,s}^\text{LiFi} = \underbrace{ (\mathbf{h}_{k,s}^\text{LiFi})^\top \mathbf{f}_{k,s}^\text{LiFi} a_{k,s}^\text{LiFi}}_{\text{Desired signal}} + \underbrace{\sum_{j \ne k} \left( \mathbf{h}_{k,s}^\text{LiFi} \right)^T \mathbf{f}_{j,s}^\text{LiFi} a_{j,s}^\text{LiFi}}_{\text{Interference from other users}} + n_{k,s}^\text{LiFi}.
		\label{eq:11}
	\end{equation}
	$n_{k,s}^\text{LiFi}$ represents the noise at the $k^\text{th}$ user at slice $s$. This noise is modeled as additive white Gaussian noise (AWGN) with zero mean and variance $(\sigma_{k,s}^\text{LiFi})^2$. Specifically, the noise for each user satisfies the Gaussian distribution $	n_{k,s}^\text{LiFi} \sim \mathcal{N}\left(0, (\sigma_{k,s}^\text{LiFi})^2 \right)$. The variance $(\sigma_{k,s}^\text{LiFi})^2$  is given by $
	(\sigma_{k,s}^\text{LiFi})^2 = \mu_{k,s}^\text{LiFi} \cdot \text{BW}_{k,s}^\text{LiFi} $, where $\mu_{k,s}^\text{LiFi}$  is the power spectral density (PSD) of AWGN for LiFi, and $\text{BW}_{k,s}^\text{LiFi}$ is the bandwidth of the LiFi link. 
	
	Furthermore, we assume that the WiFi transmitter uses a standard block-level precoding with  $\mathbf{a}_{s}^\text{WiFi} = \left[ a_{1,s}^\text{WiFi}, a_{2,s}^\text{WiFi}, \ldots, a_{K,s}^\text{WiFi} \right]^T \in \mathbb{C}^{K \times 1}$ being the symbol vector of all K users from WiFi transmitter. The symbol vector is precoded by $\mathbf{F}_{s}^\text{WiFi} = \left[ \mathbf{f}_{1,s}^\text{WiFi}, \mathbf{f}_{2,s}^\text{WiFi}, \ldots, \mathbf{f}_{K,s}^\text{WiFi} \right] \in \mathbb{C}^{M \times K}$, with $\mathbf{f}_{k,s}^\text{WiFi} = \left[ f_{1,k,s}^\text{WiFi}, f_{2,k,s}^\text{WiFi}, \ldots, f_{M,k,s}^\text{WiFi} \right]^T \in \mathbb{C}^{M \times 1}$ being the WiFi precoder .  Consequently, the transmitted signal can be written as
	\begin{equation}
		\mathbf{x}_{s}^\text{WiFi} = \mathbf{F}_{s}^\text{WiFi} \mathbf{a}_{s}^\text{WiFi} = \sum_{k=1}^{K} \mathbf{f}_{k,s}^\text{WiFi} a_{k,s}^\text{WiFi}.
		\label{eq:12}
	\end{equation}
	Thus, the signal received from the WiFi transmitter by the $k^{th}$ user at the slice $s$ is given by 
	\begin{equation}
		y_{k,s}^\text{WiFi} = \underbrace{\left( \mathbf{h}_{k,s}^\text{WiFi} \right)^\text{H} \mathbf{f}_{k,s}^\text{WiFi} a_{k,s}^\text{WiFi}}_{\text{Desired signal}} + \underbrace{\sum_{j \ne k} \left( \mathbf{h}_{k,s}^\text{WiFi} \right)^\text{H} \mathbf{f}_{j,s}^\text{WiFi} a_{j,s}^\text{WiFi}}_{\text{Interference from other users}} + n_{k,s}^\text{WiFi},
		\label{eq:13}
	\end{equation}
	where, $n_{k,s}^{WiFi}$  denotes the noise at the $k^{th}$ user at the slice $s$. The noise at each user satisfies complex complex Gaussion distribution $n_{k,s}^\text{WiFi} \sim \mathcal{CN}\left(0, (\sigma_{k,s}^\text{WiFi})^2\right)$. The variance $(\sigma_{k,s}^\text{WiFi})^2$  is given by $
	(\sigma_{k,s}^\text{WiFi})^2 = \mu_{k,s}^\text{WiFi} \cdot \text{BW}_{k,s}^\text{WiFi}$, where $\mu_{k,s}^\text{WiFi}$  is the power spectral density (PSD) of AWGN for WiFi, and $\text{BW}_{k,s}^\text{WiFi}$ is the bandwidth of the WiFi link. Each WiFi signal for the $k^{th}$ user follows a complex normal distribution with mean 0 and variance 1, denoted as $a_{k,s}^\text{WiFi} \sim \mathcal{CN}\left(0,1\right)$. In considered hybrid network, each of the \( K \) users is assigned to a specific slice based on the initial prediction stage. The interference term in the received signal equations (11) and (13) accounts for the combined effect of signals from other users on the desired signal for user \( k \). This interference can originate from users within the same slice as user \( k \) or from users in different slices due to the implementation of soft isolation among slices \cite{Soft Isolation,Dynamic Resource Sharing}. In a soft isolation scenario, where physical resources are dynamically shared among slices, there is potential for cross-slice interference. Consequently, users from different slices that share overlapping physical resources can interfere with one another, affecting the QoS experienced by user \( k \) \cite{Soft Isolation}. This interference term is therefore critical for assessing and mitigating the performance degradation caused by both intra-slice and cross-slice interference. Furthermore, the SINRs of user $k$ at the slice $s$ for WiFi and LiFi are represented as
	   %Due to the scarcity of resources, dynamic sharing among slices using soft isolation techniques is assumed \cite{Soft Isolation,Dynamic Resource Sharing}. While this approach optimizes resource utilization, it can lead to inter-slice interference when different slices simultaneously use overlapping physical resources, as well as intra-slice interference among users within the same slice \cite{Soft Isolation}. To manage these interference issues, precoding techniques are employed to optimize signal transmission and mitigate both inter-slice and intra-slice interference, thereby enhancing overall performance for users. 
	\begin{equation}
		\gamma_{k,s}^\text{WiFi} = \frac{X_{k,s}^\text{WiFi}}{Y_{k,s}^\text{WiFi}} = \frac{ \left| \left(\mathbf{h}_{k,s}^\text{WiFi} \right)^H \mathbf{f}_{k,s}^\text{WiFi} \right|^2}{\sum_{j \ne k} \left| \left(\mathbf{h}_{k,s}^\text{WiFi} \right)^H \mathbf{f}_{j,s}^\text{WiFi} \right|^2 + (\sigma_{k,s}^\text{WiFi})^2},
		\label{eq:14}
	\end{equation}
	\begin{equation}
		\gamma_{k,s}^\text{LiFi} = \frac{X_{k,s}^\text{LiFi}}{Y_{k,s}^\text{LiFi}} = \frac{ \left| \left(\mathbf{h}_{k,s}^\text{LiFi} \right)^\top \mathbf{f}_{k,s}^\text{LiFi} \right|^2}{\sum_{j \ne k} \left| \left(\mathbf{h}_{k,s}^\text{LiFi} \right)^\top \mathbf{f}_{j,s}^\text{LiFi} \right|^2 + (\sigma_{k,s}^\text{LiFi})^2}.
		\label{eq:15}
	\end{equation}
	For computing the achievable rates of LiFi links, the visible light signals must be real-valued and non-negative, making the classic Shannon formula with Gaussian input unsuitable. Since the channel capacity of LiFi is not available in closed form, we use some tight bounds to approximate the data rates. In this paper, we consider a lower bound of the achievable rate, given by \cite{LiFi Capacity}
	\begin{equation}
		R_{k,s}^\text{LiFi} = \underbrace{\frac{1}{2} \text{BW}_{k,s}^\text{LiFi} \ln \left(1 + \frac{e}{2\pi} \gamma_{k,s}^\text{LiFi}\right)}_{\text{S}_{k,s}^\text{LiFi}} \text{ nats/sec},
		\label{eq:16}
	\end{equation}
	while for WiFi, the achievable rate of user $k$ at the slice $s$ is
	\begin{equation}
		R_{k,s}^\text{WiFi} = \underbrace{\text{BW}_{k,s}^\text{WiFi} \ln \left(1 + \gamma_{k,s}^\text{WiFi}\right)}_{\text{S}_{k,s}^\text{WiFi}} \text{ nats/sec}.
		\label{eq:17}
	\end{equation}
	In our analysis, we consider three types of slices: $ s \in \{eMBB, URLLC, mMTC\}$. The equations above are only specifically applicable for $s \in eMBB$ applications, where longer block lengths are necessary to send data to meet high capacity demands. Conversely, $s \in URLLC$ requires very low latency, which can be achieved through shorter block lengths data transmission to minimize delay. For $s \in mMTC,$ which requires low latency and lower payload, short block lengths data transmissions are also preferred. Shannon capacity is not applicable for short block lengths data transmission due to the channel dispersion that arises \cite{O-RAN Network slicing}. Therefore, for the slice type of $s \in \{URLLC,mMTC\}$, which uses short block length for data transmission, the achievable rate of user $k$ at the slice $s$ using LiFi system $R_{k,s}^{LiFi}$ can be presented as 
	\begin{equation}
		R_{k,s}^\text{LiFi} = \frac{1}{2} \text{BW}_{k,s}^\text{LiFi} \left\{ 
		\underbrace{\ln \left(1 + \frac{e}{2\pi} \gamma_{k,s}^\text{LiFi}\right)}_{\text{S}_{k,s}^\text{LiFi}} 
		- \underbrace{Q^{-1} (\epsilon) \sqrt{\frac{0.5D_{k,s}^\text{LiFi}}{L_k}}}_{\text{V}_{k,s}^\text{LiFi}} 
		\right\} \text{ nats/sec},
		\label{eq:18}
	\end{equation}
	where, $ D_{k,s}^{\text{LiFi}} = 1 - \frac{1}{\left(1 + \frac{e}{2\pi} \gamma_{k,s}^{\text{LiFi}}\right)^2}$ is the channel dispersion while $L_k= BW^{LiFi}_{k,s} \times T^t_{k,s}$, indicates the associated transmit blocklength. The scaling factor of 0.5 is applied to the achievable rate and channel dispersion $D_{k,s}^{\text{LiFi}}$ due to the utilization of Hermitian symmetry \cite{VLC Short packets}. ${\text{S}_{k,s}^\text{LiFi}}$ represents the shanon term, while the ${\text{V}_{k,s}^\text{LiFi}}$ denotes the channel dispersion term. $Q^{-1}(\cdot)$ is the inverse function, with $ Q^{-1}(x) = \int_x^\infty \frac{1}{\sqrt{2\pi}} \exp\left(-\frac{t^2}{2}\right) \, dt $, and $\epsilon$ is the decoding error probability. Similarly, the achievable rate of user $ k$ at the slice $s$ using WiFi system $R_{k,s}^{\text{WiFi}}$ can be given as
	\begin{equation}
		R_{k,s}^{\text{WiFi}} = \text{BW}_{k,s}^{\text{WiFi}} \left\{ \underbrace{\ln\left(1 + \gamma_{k,s}^{\text{WiFi}}\right)}_{S_{k,s}^{\text{WiFi}}} - \underbrace{Q^{-1}(\epsilon) \sqrt{\frac{D_{k,s}^{\text{WiFi}}}{L_k}}}_{V_{k,s}^{\text{WiFi}}} \right\} \text{ nats/sec},
		\label{eq:19}
	\end{equation}
	where, $D_{k,s}^{\text{WiFi}} = 1 - \frac{1}{\left(1 + \gamma_{k,s}^{\text{WiFi}}\right)^2}$ is the channel dispersion while $L_k= BW^{WiFi}_{k,s} \times T^t_{k,s}$, indicates the associated transmit blocklength. In LiFi , total power consumption  can be modeled as 
	\begin{equation}
		P_{\text{LiFi}} = P_{\text{DC}} + P_{\text{AC}},
		\label{eq:20}
	\end{equation}
	where $P_{\text{DC}}$ and $P_{\text{AC}}$ denote the power of DC and AC current, respectively. $P_{\text{DC}}$ includes $ P_{\text{DC}}^{\text{LED}}$ for the illumination provided by the LEDs and $P_{\text{DC}}^{\text{circuit}}$ for the other circuit components. While $P_{\text{DC}}^{\text{circuit}}$ remains constant, the power used by the LEDs can be adjusted based on the dimming requirements. However, when the illumination level remains constant, $P_{\text{DC}}^{\text{LED}} = \sum_{l=1}^L V_F \cdot x_{\text{DC},l}^s$ can be considered fixed, where $V_F$ is the forward voltage of LEDs. Therefore, $P_{\text{DC}} = P_{\text{DC}}^{\text{LED}} + P_{\text{DC}}^{\text{circuit}}$ has a fixed value based on dimming requirements. The AC currents are derived directly from the output current/voltage of the precoders in LED drivers, and thus can be computed as follows \cite{AC power LiFi}
	\begin{equation}
		P_{\text{AC}} = \frac{1}{\eta^{\text{LiFi}}} \sum_{k=1}^K \left\| \boldsymbol{f}_{k,s}^{\text{LiFi}} \right\|^2,
		\label{eq:21}
	\end{equation}
	where, $\eta^{\text{LiFi}}$ signifies the power amplifier efficiency at the LiFi. In this paper, we only considers the AC power consumption for wireless signal transmission. That is, $ P_{\text{LiFi}} = P_{\text{AC}}$. For WiFi, the total power consumption can be computed as 
	\begin{equation}
		P_{\text{WiFi}} = \frac{1}{\eta^{\text{WiFi}}} \sum_{k=1}^K \left\| \boldsymbol{f}_{k,s}^{\text{WiFi}} \right\|^2 + P_c^{\text{WiFi}}.
		\label{eq:22}
	\end{equation}
	$ \eta^{\text{WiFi}}$ is the power amplifier efficiency at the WiFi. $ P_c^{\text{WiFi}}$ is the total circuit power at the WiFi. Similarly, for WiFi we considers only the power consumption for wireless signal transmission, That is $P_{\text{WiFi}} = \frac{1}{\eta^{\text{WiFi}}} \sum_{k=1}^K \left\| \boldsymbol{f}_{k,s}^{\text{WiFi}} \right\|^2$. In this work, we aim to maximize the EE of hybrid WiFi/LiFi transmitter while considering  different slice types in order to operate different applications. EE is defined as the ratio of the sum-rate of the system to the total power consumption \cite{EE Definition} as follow
	\begin{equation}
		\text{EE}^{\text{Hybrid}} = \frac{\sum_{k=1}^K R_{k,s}^{\text{WiFi}}+\sum_{k=1}^K R_{k,s}^{\text{LiFi}}}{P_{\text{WiFi}}+P_{\text{LiFi}}}.
		\label{eq:23}
	\end{equation}
	%For URLLC and mMTC services, we assume that the \( k \)-th IIoT device or transmitter follows an independent and identically distributed Poisson distribution with a packet arrival rate \( \lambda \) and a data packet size \( L_{\text{packet}} \) in bytes \cite{27}.
	Moreover, the total latency \( T_{k,s}\) for transmission to  the $k^{th}$ user at the slice $s$ includes several components: the queue waiting time \( T^w_{k,s} \), the transmission time \( T^t_{k,s} \), the channel access delay \( T^a_{k,s} \), the backhaul delay \( T^b_{k,s} \), the packet reception delay \( T^r_{k,s} \), and the processing delay \( T^p_{k,s} \). These components together define the total latency, given by \cite{Latency} 
	
	\begin{equation}
	T_{k,s} = T^w_{k,s} + T^t_{k,s} + T^a_{k,s} + T^b_{k,s} + T^r_{k,s} + T^p_{k,s}.
	\label{eq:24}
    \end{equation}
	
	\section{THE SOLUTION OF FORMULATED PROBLEM}
	The main objective is to maximize the EE at Hybrid WiFi/LiFi transmitter for any slice type $s$ predicted for the users in the stage 1 in order to operate different applications. In this stage 2, the EE of the system is maximized through the optimization of the transmit precoding under required system constraints. The Hybrid system constraints include the constraint on LEDs input signals for LiFi, end-to-end latency constraint,  data rate constraint, and power budget constraint for Hybrid system. Hence, the EE optimization problem is formulated as
	\begin{equation}
		\text{P1:} \quad \max_{\boldsymbol{f}^{\text{WiFi}},\boldsymbol{f}^{\text{LiFi}}} \text{EE}^{\text{Hybrid}} = \frac{\sum_{k=1}^K R_{k,s}^{\text{WiFi}}+\sum_{k=1}^K R_{k,s}^{\text{LiFi}}}{P_{\text{WiFi}}+P_{\text{LiFi}}} 
		\label{eq:25}
	\end{equation}
	
	Subject to:
	
	\begin{equation}
		\begin{aligned}
			& \text{C1:} \quad R_{k,s}^{\text{WiFi}} + R_{k,s}^{\text{LiFi}} \geq R_{k,s}^{\text{min}}, \quad \forall k, \\
			& \text{C2:} \quad P_{\text{WiFi}} + P_{\text{LiFi}} \leq P_{\text{max}}^{\text{Hybrid}},\\
			& \text{C3:} \quad  \sum_{k=1}^K \left| f_{l,k,s}^{\text{LiFi}} \right| \leq \min \left( x_{DC,l}^s, I_{\text{max},s} - x_{DC,l}^s \right), \quad \forall l,\\
			& \text{C4:} \quad T_{k,s} \leq T^{min}_{k,s}, \quad \forall k.
			\label{eq:26}	
		\end{aligned}
	\end{equation}
	The optimization problem $\text{P1}$ is non-convex, making it challenging to achieve a globally optimal solution in practice. To address this, a sequential convex programming (SCP) strategy is employed to obtain a provably convergent iterative convex solution. By using SCP, the problem P1 are reformulated into a more manageable form by introducing scalar variable: $\psi^\text{Hybrid}$, representing the total power consumption of the hybrid WiFi/LiFi system, and $\phi^{\text{Hybrid}}$ , representing the EE of the hybrid system. The equivalent transformed optimization problem for P1 can then be expressed in P2 as follows
	\begin{equation}
		\text{P2:} \quad \max_{\boldsymbol{f}^{\text{WiFi}},\boldsymbol{f}^{\text{LiFi}}, \phi^{\text{Hybrid}}, \psi^{\text{Hybrid}}} \; \phi^{\text{Hybrid}}
		\label{eq:27}
	\end{equation}
	
	Subject to:
	
	\begin{equation}
		\begin{aligned}
			& \text{C5:} \quad \text{EE}^{\text{Hybrid}} \geq \phi^{\text{Hybrid}}, \\
			& \text{C6:} \quad \sum_{k=1}^K R_{k,s}^{\text{WiFi}}+\sum_{k=1}^K R_{k,s}^{\text{LiFi}} \geq \phi^{\text{Hybrid}} \psi^{\text{Hybrid}}, \\
			& \text{C7:} \quad P_{\text{WiFi}} + P_{\text{LiFi}} \leq \psi^{\text{Hybrid}},\\
			&  \quad \quad \: \: \text{C1}, \text{C2}, \text{C3} \quad \text{and}\quad \text{C4}.
			\label{eq:28}
		\end{aligned}
	\end{equation}
	In the optimization problem P2, constraints C6 and C7 are introduced to handle the non-convex constraint C5. Constraint C6 is non-convex because of the coupled optimization variables $\phi^{\text{Hybrid}} \psi^\text{Hybrid}$ and its dependence on the rate expressions $R_{k,s}^{\text{WiFi}} $ and  $R_{k,s}^{\text{LiFi}}$. We exploited the first-order Taylor expansion to approximate the convexity of constraint C6 in problem P2 in terms of the coupled optimization variables $\phi^{\text{Hybrid}} \psi^\text{Hybrid}$ as
	\begin{equation}
		\sum_{k=1}^K R_{k,s}^{\text{WiFi}}+\sum_{k=1}^K R_{k,s}^{\text{LiFi}} \geq \phi^\text{Hybrid}_{(t)} \psi^\text{Hybrid}_{(t)}
		+ \psi^\text{Hybrid}_{(t)} \left( \phi^\text{Hybrid} - \phi^\text{Hybrid}_{(t)} \right) 
		+ \phi^\text{Hybrid}_{(t)} \left( \psi^\text{Hybrid} - \psi^\text{Hybrid}_{(t)} \right),
		\label{eq:29}
	\end{equation}
	where $\phi^\text{Hybrid}_{(t)} \psi^\text{Hybrid}_{(t)}$ indicate the value of $\phi^\text{Hybrid} \psi^\text{Hybrid}$ in the $t^{th}$ iteration. Further, to deal with the non-convexity of $\text{C6}$ in $\text{P2}$, it is important to approximate the non-convex rate expressions $R_{k,s}^\text{WiFi}$ and $R_{k,s}^\text{LiFi}$. For \( s \in {eMBB} \), the rate expressions have only the Shannon term, that is  $S_{k,s}^\text{WiFi}$ and $S_{k,s}^\text{LiFi}$. While for \( s \in \{ {URLLC}, {mMTC} \} \), the rate expressions have both the Shannon term and the channel dispersion term, that is $ S_{k,s}^\text{WiFi}$, $V_{k,s}^\text{WiFi}$, $S_{k,s}^\text{LiFi}$, and $V_{k,s}^\text{LiFi}$. To convexify \( \text{C6} \) in terms of rate expressions, it is vital to obtain a concave approximation of Shannon’s term and a convex approximation of the channel dispersion term.
	\subsection{Concave approximation for $S_{k,s}^\text{WiFi}$ and $S_{k,s}^\text{LiFi}$ }
	
	We first obtained the concave approximation for $S_{k,s}^{WiFi}$ and $S_{k,s}^{LiFi}$ as function of $\boldsymbol{f}_{k,s}^{WiFi}$ and $\boldsymbol{f}_{k,s}^{LiFi}$. According to [\cite{Shanon term approximation}, Eq. (20)], a concave approximation of $S_{k,s}^{WiFi}$ and $S_{k,s}^{LiFi}$ can be given as 
	\begin{align}
		S_{k,s}^\text{WiFi} &\geq \left[ S_{k,s}^\text{WiFi} \right]^{(t)} - \left[ \gamma_{k,s}^\text{WiFi} \right]^{(t)} 
		+ A_{k,s}^{WiFi}-B_{k,s}^{WiFi}, \nonumber \\
		&\quad  \overset{\Delta}{=} \overline {S_{k,s}^\text{WiFi}},
		\label{eq:30}
	\end{align}
	\begin{align}
		S_{k,s}^\text{LiFi} & \geq \left[ S_{k,s}^\text{LiFi} \right]^{(t)} - \frac{e}{2\pi} \left[ \gamma_{k,s}^\text{LiFi} \right]^{(t)} + A_{k,s}^{LiFi}-B_{k,s}^{LiFi} , \nonumber \\
		&  \overset{\Delta}{=} \overline {S_{k,s}^\text{LiFi}},
		\label{eq:31}
	\end{align}
	where 
	\begin{align}
		&A_{k,s}^{WiFi}=\frac{2 \Re \left\{ \left(  \left(\mathbf{h}_{k,s}^\text{WiFi} \right)^H (\mathbf{f}_{k,s}^\text{WiFi})^{(t)}  \right)^H \left(\left(\mathbf{h}_{k,s}^\text{WiFi} \right)^H (\mathbf{f}_{k,s}^\text{WiFi})\right) \right\}}{\left[ Y_{k,s}^\text{WiFi} \right]^{(t)}}, \nonumber\\
		&B_{k,s}^{WiFi}=\frac{ \left[ X_{k,s}^\text{WiFi} \right]^{(t)} }{\left[ Y_{k,s}^\text{WiFi} \right]^{(t)} \left( \left[ Y_{k,s}^\text{WiFi} \right]^{(t)} +  \left[ X_{k,s}^\text{WiFi} \right]^{(t)}  \right)} \left( \left[ X_{k,s}^\text{WiFi} \right] + \left[ Y_{k,s}^\text{WiFi} \right] \right),
	\end{align}
	\begin{align}
		&A_{k,s}^{LiFi}=\frac{2e \Re \left\{ \left(  \left(\mathbf{h}_{k,s}^\text{LiFi} \right)^\top (\mathbf{f}_{k,s}^\text{LiFi})^{(t)}  \right)^\top \left(\left(\mathbf{h}_{k,s}^\text{LiFi} \right)^\top (\mathbf{f}_{k,s}^\text{LiFi})\right) \right\}}{2\pi \left[ Y_{k,s}^\text{LiFi} \right]^{(t)}},  \nonumber\\
		&B_{k,s}^{LiFi}=\frac{ e \left[ X_{k,s}^\text{LiFi} \right]^{(t)} }{\left( 2\pi \left[ Y_{k,s}^\text{LiFi} \right]^{(t)} \right) \left( 2\pi \left[ Y_{k,s}^\text{LiFi} \right]^{(t)} +  e \left[ X_{k,s}^\text{LiFi} \right]^{(t)}  \right)} \left(  \left[e X_{k,s}^\text{LiFi} \right] + 2\pi \left[ Y_{k,s}^\text{LiFi} \right] \right).
	\end{align}
	$[\gamma_{k,s}^\text{WiFi}]^{(t)}$, $[\gamma_{k,s}^\text{LiFi}]^{ \left(t\right)}$, $[S_{k,s}^\text{WiFi}]^{(t)}$,   $ [S_{k,s}^\text{LiFi}]^{\left(t\right)}$, $\left[ X_{k,s}^\text{WiFi} \right]^{(t)}$,$\left[ X_{k,s}^\text{LiFi} \right]^{(t)}$,$\left[ Y_{k,s}^\text{WiFi} \right]^{(t)}$, and $\left[ Y_{k,s}^\text{LiFi} \right]^{(t)}$ indicate the values of $\gamma_{k,s}^\text{WiFi}$, $\gamma_{k,s}^{LiFi}$, $S_{k,s}^{WiFi}$,$S_{k,s}^{LiFi}$, $X_{k,s}^\text{WiFi}$, $X_{k,s}^\text{LiFi}$, $Y_{k,s}^\text{WiFi}$, and $Y_{k,s}^\text{LiFi}$ in the $t^{th}$ iteration, respectively. It is under the trust region constrained by
	\begin{equation}
		\left| (\mathbf{h}_{k,s}^\text{WiFi}) ^H \mathbf{f}_{k,s}^\text{WiFi} \right|^2 > 0,
	\end{equation}
	\begin{equation}
		e*\left| (\mathbf{h}_{k,s}^\text{LiFi})^\top \mathbf{f}_{k,s}^\text{LiFi} \right|^2 > 0,
	\end{equation}
	where $\left| (\mathbf{h}_{k,s}^\text{WiFi}) ^H \mathbf{f}_{k,s}^\text{WiFi} \right|^2$ and $\left| (\mathbf{h}_{k,s}^\text{LiFi})^\top \mathbf{f}_{k,s}^\text{LiFi} \right|^2$ can be approximated by the first order taylor approximation as
	\begin{align}
		&2 \Re \left[ \left( \left( \mathbf{h}_{k,s}^\text{WiFi} \right)^H (\mathbf{f}_{k,s}^\text{WiFi})^{(t)} \right)^H \left( \mathbf{h}_{k,s}^\text{WiFi} \right)^H \mathbf{f}_{k,s}^\text{WiFi} \right]- \left| \left( \mathbf{h}_{k,s}^\text{WiFi} \right)^H \mathbf{f}_{k,s}^\text{WiFi} \right|^2 > 0,\\
		& e\Big\{2 \Re \left[ \left( \left( \mathbf{h}_{k,s}^\text{LiFi} \right)^\top (\mathbf{f}_{k,s}^\text{LiFi})^{(t)} \right)^\top \left( \mathbf{h}_{k,s}^\text{LiFi} \right)^\top \mathbf{f}_{k,s}^\text{LiFi} \right]- \left| \left( \mathbf{h}_{k,s}^\text{LiFi} \right)^\top \mathbf{f}_{k,s}^\text{LiFi} \right|^2 \Big\} > 0 .
	\end{align}
	
	\subsection{Convex approximation for $V_{k,s}^\text{WiFi}$ and $V_{k,s}^\text{LiFi}$ }
	We can write $V_{k,s}^\text{WiFi}$ and $V_{k,s}^\text{LiFi}$ as
	\begin{equation}
		V_{k,s}^\text{WiFi} = \chi_{k,s}^\text{WiFi} \sqrt{D_{k,s}^\text{WiFi}},
	\end{equation}
	and
	\begin{equation}
		V_{k,s}^\text{LiFi} = \chi_{k,s}^\text{LiFi} \sqrt{0.5  D_{(k,s)}^\text{LiFi}},	
	\end{equation}
	where, $\chi_{k,s}^\text{WiFi}$ and $\chi_{k,s}^\text{LiFi} = \frac{Q^{-1}(\epsilon)}{\sqrt{L_k}}$. $ D_{k,s}^\text{WiFi} = 1 - \frac{1}{(1 + \gamma_{k,s}^\text{WiFi})^2}$, and $ D_{k,s}^\text{LiFi} = 1 - \frac{1}{\left(1 + \frac{e}{2\pi} \gamma_{k,s}^\text{LiFi}\right)^2}$. In Eq. (38) and (39), the optimization variables exist in terms of $ D_{k,s}^\text{WiFi}$ and $D_{k,s}^\text{LiFi}$. Both terms are in the form $ \sqrt{x}$. We found that $f(x) = \sqrt{x}$ is concave for $ x > 0$, because its second derivative is negative for all $ x > 0$. Therefore, its convex approximation around a feasible point $ x^{(t)}$ can be computed using its first-order Taylor approximation as
	\begin{equation}
		\sqrt{x} \leq \sqrt{x^{(t)}} + \frac{1}{2 \sqrt{x^{(t)}}} \left( x - x^{(t)} \right),
	\end{equation}
	\begin{equation}
		\sqrt{x} \leq \frac{\sqrt{x^{(t)}}}{2} + \frac{x}{2 \sqrt{x^{(t)}}}.
	\end{equation}
	The above inequality holds true for all $x > 0$ and $ x^{(t)} > 0 $. By using Eq. (41) we have
	\begin{equation}
		\sqrt{D_{k,s}^\text{WiFi}} \leq \frac{\sqrt{D_{k,s}^\text{WiFi(t)}}}{2} + \frac{D_{k,s}^\text{WiFi}}{2 \sqrt{D_{k,s}^\text{WiFi(t)}}},
	\end{equation}
	\begin{equation}
		\sqrt{D_{k,s}^\text{WiFi}} \leq \frac{\sqrt{D_{k,s}^\text{WiFi(t)}}}{2} + \frac{1}{2 \sqrt{D_{k,s}^\text{WiFi(t)}}} - \frac{1}{2 \sqrt{D_{k,s}^\text{WiFi(t)}}} \left( \frac{1}{(1 + \gamma_{k,s}^\text{WiFi})^2} \right),
	\end{equation} 
	\begin{equation}
		\sqrt{D_{k,s}^\text{WiFi}} \leq \frac{\sqrt{D_{k,s}^\text{WiFi(t)}}}{2} + \frac{1}{2 \sqrt{D_{k,s}^\text{WiFi(t)}}} - \frac{1}{2 \sqrt{D_{k,s}^\text{RF(t)}}} \Bigg[\frac{\Big(\sum_{j \ne k} \left| \left(\mathbf{h}_{k,s}^\text{WiFi} \right)^H \mathbf{f}_{j,s}^\text{WiFi} \right|^2 + (\sigma_{k,s}^\text{WiFi})^2\Big)^2}{\Big(\sum_{k=1}^K \left| \left(\mathbf{h}_{k,s}^\text{WiFi} \right)^H \mathbf{f}_{k,s}^\text{WiFi} \right|^2 + (\sigma_{k,s}^\text{WiFi})^2\Big)^2}\Bigg],
	\end{equation}
	
	\begin{equation}
		\sqrt{D_{k,s}^\text{WiFi}} \leq \frac{\sqrt{D_{k,s}^\text{WiFi(t)}}}{2} + \frac{1}{2 \sqrt{D_{k,s}^\text{WiFi(t)}}} - \frac{1}{2 \sqrt{D_{k,s}^\text{WiFi(t)}}} \left( \frac{(Y_{k,s}^\text{WiFi})^2}{(X_{k,s}^\text{WiFi} + Y_{k,s}^\text{WiFi})^2} \right),
	\end{equation}
	\begin{equation}
		\sqrt{D_{k,s}^\text{WiFi}} \leq \alpha_{k,s}^\text{WiFi (t)} - \beta_{k,s}^\text{WiFi (t)} \left( \frac{(Y_{k,s}^\text{WiFi})^2}{(X_{k,s}^\text{WiFi} + Y_{k,s}^\text{WiFi})^2} \right),
	\end{equation}
	for
	\begin{equation}
		0 < \beta_{k,s}^\text{WiFi (t)} = \frac{1}{2 \sqrt{D_{k,s}^\text{WiFi (t)}}},
	\end{equation}
	\begin{equation}
		0 < \alpha_{k,s}^\text{WiFi (t)} = \frac{\sqrt{D_{k,s}^\text{WiFi (t)}}}{2} + \beta_{k,s}^\text{WiFi (t)}.
	\end{equation}
	The right-hand side (R.H.S) of Eq. (46) is still non-convex due to the term $\frac{(Y_{k,s}^\text{WiFi})^2}{(X_{k,s}^\text{WiFi} + Y_{k,s}^\text{WiFi})^2}$
	therefore, we obtain the following concave approximation for it
	\begin{equation}
		\frac{(Y_{k,s}^\text{WiFi})^2}{(X_{k,s}^\text{WiFi} + Y_{k,s}^\text{WiFi})^2} = \frac{(Y_{k,s}^\text{WiFi})^2}{(X_{k,s}^\text{WiFi} + Y_{k,s}^\text{WiFi})} \times \frac{1}{(X_{k,s}^\text{WiFi} + Y_{k,s}^\text{WiFi})}.
	\end{equation}
	We found that $ f(x) = \frac{1}{x}$ is convex for $x > 0$, because its second derivative is positive for all $ x > 0$. The following inequality holds true for all $ x > 0 $ and $ x^{(t)} > 0 $
	\begin{equation}
		\frac{1}{x} \geq \frac{1}{x^{(t)}} + \frac{x - x^{(t)}}{(x^{(t)})^2},
	\end{equation}
	\begin{equation}
		\frac{1}{x} \geq \frac{2}{x^{(t)}} - \frac{x}{(x^{(t)})^2}.
	\end{equation}
	Then, 
	\begin{equation}
		\frac{1}{X_{k,s}^{\text{WiFi}} + Y_{k,s}^{\text{WiFi}}} = \frac{2}{X_{k,s}^{WiFi(t)} + Y_{k,s}^{WiFi(t)}} - \frac{X_{k,s}^{\text{WiFi}} + Y_{k,s}^{\text{WiFi}}}{(X_{k,s}^{WiFi(t)} + Y_{k,s}^{WiFi(t)})^2},
	\end{equation}
	Under the condition 
	\begin{equation}
		X_{k,s}^{\text{WiFi}} + Y_{k,s}^{\text{WiFi}} \leq 2 \left( X_{k,s}^{WiFi(t)} + Y_{k,s}^{WiFi(t)} \right).
	\end{equation}
	Putting Eq. (52) in Eq. (49) results in
	\begin{equation}
		\frac{(Y_{k,s}^{\text{WiFi}})^2}{(X_{k,s}^{\text{WiFi}} + Y_{k,s}^{\text{WiFi}})^2} = \frac{(Y_{k,s}^{\text{WiFi}})^2}{(X_{k,s}^{\text{WiFi}} + Y_{k,s}^{\text{WiFi}})} \times \left( \frac{2}{X_{k,s}^{WiFi(t)} + Y_{k,s}^{WiFi(t)}} - \frac{X_{k,s}^{\text{WiFi}} + Y_{k,s}^{\text{WiFi}}}{(X_{k,s}^{WiFi(t)} + Y_{k,s}^{WiFi(t)})^2} \right),
	\end{equation}
	\begin{equation}
		=\frac{(Y_{k,s}^{\text{WiFi}})^2}{X_{k,s}^{\text{WiFi}} + Y_{k,s}^{\text{WiFi}}} \times \frac{2}{X_{k,s}^{WiFi(t)} + Y_{k,s}^{WiFi(t)}} - \frac{(Y_{k,s}^{\text{WiFi}})^2}{(X_{k,s}^{WiFi(t)} + Y_{k,s}^{WiFi(t)})^2}.
	\end{equation}
	We now employ the inequality \text{[\cite{Shanon term approximation}, Eq. 21]} to approximate
	$\frac{(Y_{k,s}^{\text{WiFi}})^2}{X_{k,s}^{\text{WiFi}} + Y_{k,s}^{\text{WiFi}}}$
	\begin{equation}
		\frac{(Y_{k,s}^{\text{WiFi}})^2}{X_{k,s}^{\text{WiFi}} + Y_{k,s}^{\text{WiFi}}} \geq \frac{2 Y_{k,s}^{\text{WiFi(t)}}}{X_{k,s}^{\text{WiFi(t)}} + Y_{k,s}^{\text{WiFi(t)}}} \cdot Y_{k,s}^{\text{WiFi}} - \frac{(Y_{k,s}^{\text{WiFi}})^{(t)})^2}{(X_{k,s}^{\text{WiFi(t)}} + Y_{k,s}^{\text{WiFi(t)}})^{2}} \cdot (X_{k,s}^{\text{WiFi}} + Y_{k,s}^{\text{WiFi}}),
	\end{equation}
	\begin{equation}
		\frac{(Y_{k,s}^{\text{WiFi}})^2}{X_{k,s}^{\text{WiFi}} + Y_{k,s}^{\text{WiFi}}} \geq \frac{2 Y_{k,s}^{\text{WiFi(t)}}}{X_{k,s}^{\text{WiFi(t)}} + Y_{k,s}^{\text{WiFi(t)}}} \cdot \tau_{k,s}^{\text{WiFi(t)}} - \frac{(Y_{k,s}^{\text{WiFi}})^{(t)})^2}{(X_{k,s}^{\text{WiFi(t)}} + Y_{k,s}^{\text{WiFi(t)}})^2} \cdot (X_{k,s}^{\text{WiFi}} + Y_{k,s}^{\text{WiFi}}),
	\end{equation}
	Under the condition 
	\begin{equation}
		X_{k,s}^{\text{WiFi}} + Y_{k,s}^{\text{WiFi}} \leq 2 \cdot \frac{X_{k,s}^{\text{WiFi(t)}} + Y_{k,s}^{\text{WiFi(t)}}}{Y_{k,s}^{\text{WiFi(t)}}} \cdot Y_{k,s}^{\text{WiFi}},
	\end{equation}
	where $\tau_{k,s}^{\text{WiFi}(t)}$ is the first-order Taylor approximation of $Y_{k,s}^{\text{WiFi}}$
	\begin{equation}
		Y_{k,s}^{\text{WiFi}} \geq \tau_{k,s}^{\text{WiFi(t)}} = \sum_{j \ne k} \left( 2 \Re \left\{\Big( \left( \mathbf{h}_{k,s}^{\text{WiFi}} \right)^H \mathbf{f}_{j,s}^{\text{WiFi}(t)}\Big)^* \left( \mathbf{h}_{k,s}^{\text{WiFi}} \right)^H \mathbf{f}_{k,s}^{\text{WiFi}}\right\} - \left| \left( \mathbf{h}_{k,s}^{\text{WiFi}} \right)^H \mathbf{f}_{j,s}^{\text{WiFi}} \right|^2 \right) + \sigma_{k,s}^{\text{WiFi2}}.
	\end{equation}
	Putting Eq. (57) in Eq. (55) results in 
	\begin{align}
		\frac{\left( Y_{k,s}^{\text{WiFi}} \right)^2}{\left( X_{k,s}^{\text{WiFi}} + Y_{k,s}^{\text{WiFi}} \right)^2} &= \left( \frac{2 Y_{k,s}^{\text{WiFi}(t)}}{X_{k,s}^{\text{WiFi}(t)} + Y_{k,s}^{\text{WiFi}(t)}} \cdot \tau_{k,s}^{\text{WiFi}(t)} - \frac{\left( Y_{k,s}^{\text{WiFi}(t)} \right)^2}{\left( X_{k,s}^{\text{WiFi}(t)} + Y_{k,s}^{\text{WiFi}(t)} \right)^2} \cdot \left( X_{k,s}^{\text{WiFi}} + Y_{k,s}^{\text{WiFi}} \right) \right) \nonumber \\
		&\times \frac{2}{X_{k,s}^{\text{WiFi}(t)} + Y_{k,s}^{\text{WiFi}(t)}} - \frac{\left( Y_{k,s}^{\text{WiFi}} \right)^2}{\left( X_{k,s}^{\text{WiFi}(t)} + Y_{k,s}^{\text{WiFi}(t)} \right)^2},
	\end{align}
	\begin{align}
		\frac{\left( Y_{k,s}^{\text{WiFi}} \right)^2}{\left( X_{k,s}^{\text{WiFi}} + Y_{k,s}^{\text{WiFi}} \right)^2} &= \frac{4 Y_{k,s}^{\text{WiFi}(t)}}{\left( X_{k,s}^{\text{WiFi}(t)} + Y_{k,s}^{\text{WiFi}(t)} \right)^2} \cdot \tau_{k,s}^{\text{WiFi}(t)} - \frac{2 \left( Y_{k,s}^{\text{WiFi}(t)} \right)^2}{\left( X_{k,s}^{\text{WiFi}(t)} + Y_{k,s}^{\text{WiFi}(t)} \right)^3} \cdot \left( X_{k,s}^{\text{WiFi}} + Y_{k,s}^{\text{WiFi}} \right) \nonumber \\
		&- \frac{\left( Y_{k,s}^{\text{WiFi}} \right)^2}{\left( X_{k,s}^{\text{WiFi}(t)} + Y_{k,s}^{\text{WiFi}(t)} \right)^2}.
	\end{align}
Substituting $Eq. (61)$ into $Eq. (46)$ results in Convex approximationfor $D_{k,s}^{WiFi}$ as follows 
\begin{align}
	\sqrt{D_{k,s}^{\text{WiFi}}} &\leq \alpha_{k,s}^{\text{WiFi}(t)} - \beta_{k,s}^{\text{WiFi}(t)} \bigg[ \frac{4 Y_{k,s}^{\text{WiFi}(t)}}{\left( X_{k,s}^{\text{WiFi}(t)} + Y_{k,s}^{\text{WiFi}(t)} \right)^2} \cdot \tau_{k,s}^{\text{WiFi}(t)}- \frac{2 \left( Y_{k,s}^{\text{WiFi}(t)} \right)^2}{\left( X_{k,s}^{\text{WiFi}(t)} + Y_{k,s}^{\text{WiFi}(t)} \right)^3}\nonumber \\
	 &\times \left( X_{k,s}^{\text{WiFi}} + Y_{k,s}^{\text{WiFi}} \right)- \frac{\left( Y_{k,s}^{\text{WiFi}} \right)^2}{\left( X_{k,s}^{\text{WiFi}(t)} + Y_{k,s}^{\text{WiFi}(t)} \right)^2} \bigg], \nonumber \\
	 &  \overset{\Delta}{=}\sqrt{\overline {D_{k,s}^\text{WiFi}}}. 
\end{align}
Subsequently, 
\begin{equation}
	\overline{V_{k,s}^{\text{WiFi}}} = \chi_{k,s}^{\text{WiFi}} \sqrt{\overline{D_{k,s}^{\text{WiFi}}}}.
\end{equation}
Similarly, we derived convex approximation for $D_{k,s}^{LiFi}$ as
\begin{align}
	\sqrt{D_{k,s}^{\text{LiFi}}} &\leq \alpha_{k,s}^{\text{LiFi}(t)} - \beta_{k,s}^{\text{LiFi}(t)} \bigg[ \frac{4 Y_{k,s}^{\text{LiFi}(t)}}{\left( X_{k,s}^{\text{LiFi}(t)} + Y_{k,s}^{\text{LiFi}(t)} \right)^2} \cdot \tau_{k,s}^{\text{LiFi}(t)}- \frac{2 \left( Y_{k,s}^{\text{LiFi}(t)} \right)^2}{\left( X_{k,s}^{\text{LiFi}(t)} + Y_{k,s}^{\text{LiFi}(t)} \right)^3}\nonumber \\
	&\times \left( eX_{k,s}^{\text{LiFi}} +2\pi Y_{k,s}^{\text{LiFi}} \right)- \frac{\left( 2\pi Y_{k,s}^{\text{LiFi}} \right)^2}{\left( X_{k,s}^{\text{LiFi}(t)} + Y_{k,s}^{\text{LiFi}(t)} \right)^2} \bigg], \nonumber \\
	&  \overset{\Delta}{=}\sqrt{\overline {D_{k,s}^\text{LiFi}}}. 
\end{align}
for
	\begin{equation}
	0 < \beta_{k,s}^\text{LiFi (t)} = \frac{1}{2 \sqrt{D_{k,s}^\text{LiFi (t)}}},
\end{equation}
\begin{equation}
	0 < \alpha_{k,s}^\text{LiFi (t)} = \frac{\sqrt{D_{k,s}^\text{LiFi (t)}}}{2} + \beta_{k,s}^\text{LiFi (t)},
\end{equation}
and under the condition 
\begin{equation}
	eX_{k,s}^{\text{LiFi}} + 2\pi Y_{k,s}^{\text{LiFi}} \leq 2 \left( X_{k,s}^{LiFi(t)} + Y_{k,s}^{LiFi(t)} \right),
\end{equation}
\begin{equation}
	eX_{k,s}^{\text{LiFi}} + 2\pi Y_{k,s}^{\text{LiFi}} \leq 2 \cdot \frac{X_{k,s}^{\text{LiFi(t)}} + Y_{k,s}^{\text{LiFi(t)}}}{Y_{k,s}^{\text{LiFi(t)}}} \cdot 2\pi Y_{k,s}^{\text{LiFi}},
\end{equation}
where, 
\begin{align}
&D_{k,s}^{\text{LiFi(t)}} = 1 - \frac{1}{\left(1 + \frac{e}{2\pi} \gamma_{k,s}^{\text{LiFi(t)}}\right)^2},\\
&Y_{k,s}^{\text{LiFi}(t)}=2 \pi \Big({\sum_{j \ne k} \left| \left(\mathbf{h}_{k,s}^\text{LiFi} \right)^\top \mathbf{f}_{j,s}^\text{LiFi(t)} \right|^2 + (\sigma_{k,s}^\text{LiFi})^2}\Big),\\
&X_{k,s}^{\text{LiFi}(t)}=e\left| \left(\mathbf{h}_{k,s}^\text{LiFi} \right)^\top \mathbf{f}_{k,s}^\text{LiFi} \right|^2, \\
&\tau_{k,s}^{\text{LiFi(t)}} = 2 \pi \Big[\sum_{j \ne k} \left( 2 \Re \left\{\Big( \left( \mathbf{h}_{k,s}^{\text{LiFi}} \right)^\top \mathbf{f}_{j,s}^{\text{LiFi}(t)}\Big)^* \left( \mathbf{h}_{k,s}^{\text{LiFi}} \right)^H \mathbf{f}_{k,s}^{\text{LiFi}}\right\} - \left| \left( \mathbf{h}_{k,s}^{\text{LiFi}} \right)^H \mathbf{f}_{j,s}^{\text{LiFi}} \right|^2 \right) + \sigma_{k,s}^{\text{LiFi2}}\Big].
\end{align}
Subsequently, 
\begin{equation}
	\overline{V_{k,s}^{\text{LiFi}}} = \chi_{k,s}^{\text{LiFi}} \sqrt{\overline{D_{k,s}^{\text{LiFi}}}}.
\end{equation}
Now, constraint C6 in P2 in term of rate expressions is convex, because we obtained concave approximation of shanon’s term and convex approximation of channel dispersion term.  C6 in convex form can be represented as follows
\begin{equation}
\sum_{k=1}^K \overline{R_{k,s}^{\text{WiFi}}}+\sum_{k=1}^K \overline{R_{k,s}^{\text{LiFi}}} \geq \phi^\text{Hybrid}_{(t)} \psi^\text{Hybrid}_{(t)}
	+ \psi^\text{Hybrid}_{(t)} \left( \phi^\text{Hybrid} - \phi^\text{Hybrid}_{(t)} \right) 
	+ \phi^\text{Hybrid}_{(t)} \left( \psi^\text{Hybrid} - \psi^\text{Hybrid}_{(t)} \right),
\end{equation}
Where, \(\overline{R_{k,s}^{\text{WiFi}}} = \overline{S_{k,s}^{\text{WiFi}}} - \overline{V_{k,s}^{\text{WiFi}}}\) and \(\overline{R_{k,s}^{\text{LiFi}}} = \overline{S_{k,s}^{\text{LiFi}}} - \overline{V_{k,s}^{\text{LiFi}}}\). Furthermore, C1 in P2 is now convex and can be expressed as
\begin{equation}
	\overline{R_{k,s}^{\text{WiFi}}}+ \overline{R_{k,s}^{\text{LiFi}}} \geq R_{k,s}^{\text{min}}.
\end{equation}
Eventually, the optimization problem P2  presented in Eq. (27)  can be reformulated as following convex optimization problems P3 as
	
	\begin{equation}
		\text{P3:} \quad \max_{\boldsymbol{f}^{\text{WiFi}},\boldsymbol{f}^{\text{LiFi}}, \phi^{\text{Hybrid}}, \psi^{\text{Hybrid}}} \; \phi^{\text{Hybrid}}
	\end{equation}
	
	Subject to:
	\begin{equation}
		\begin{aligned}
			&  \quad \quad \: \:  \text{C2}, \text{C3},\text{C4}, \text{C7}, (36), (37), (47), (48),(53), (58), (65), (66), (67), (68), (74) \quad \text{and} \quad (75)  .
		\end{aligned}
	\end{equation}
	
The reformulated optimization problems in Eq. (76)  has been made convex and can now be efficiently solved using MATLAB CVX at the $t$ iteration. This iterative approach ultimately converges to a suboptimal solution for the original problems. Proper scaling of SINR terms is crucial before applying CVX optimization, as it prevents numerical issues, ensuring accurate and stable convergence to optimal solutions.  The proposed algorithm is applicable for EE maximization with data transmission from a single transmitter as well, whether WiFi or LiFi. The final sub-optimal solution of the proposed algorithm depends on the initial point, which must be feasible and can be selected either randomly or heuristically \cite{Initial Point}. However, using the finite block length rates as given in Eqs. (18) and (19) can lead to negative rates for $\gamma_{k,s}^\text{WiFi}<<1$ and $\gamma_{k,s}^\text{LiFi}<<1$. Moreover, the data rate constraint in Eq. (76) may not be satisfied with a random initial point. To ensure a feasible initial point, we utilize the approach described in Appendix.

\subsection{Computational Complexity Analysis of the Proposed Algorithm:}
The computational complexity of the proposed predict-and-optimize algorithm for hybrid WiFi/LiFi networks is analyzed in two distinct stages: prediction of network slices and optimization of precoding for EE. In the first stage, the computational complexity of using resilient backpropagation for predicting slice types (eMBB, URLLC, mMTC) in a neural network trained on historically monitored KPIs primarily depends on the forward pass through the network. For a model with \(L\) layers, where each layer \(l\) has \(n_l\) neurons, the complexity is \(O\left(m \cdot n_1 + \sum_{l=1}^{L-1} n_l \cdot n_{l+1}\right)\), with \(m\) representing the input KPIs. This reflects the operations needed to compute the output from the input features, making the model efficient for real-time predictions in hybrid WiFi/LiFi networks. In the second stage, optimization of precoding, the complexity arises from applying sequential convex approximation and inner approximation methods. The convex problem presented in  Eq. (76) involves $2KM+2KL+2$ scalar optimization variables and $12K+L+3$ linear constraints and second order cone (SOC) constraints. Consequently, the worst-case per-iteration complexity of this stage, using the interior-point method, is estimated as $O\Big(\sqrt{12K+L+3}\big(2KM+2KL+2\big)^3\Big)$.

\section{Numerical Simulations}	
In this section, we detail the simulation setup used to evaluate our proposed predict and optimize algorithm for energy-efficient service differentiation in a hybrid LiFi and WiFi network. The simulation environment models a room with dimensions of 10 meters in length, 10 meters in width, and 5 meters in height, equipped with both multi-LED LiFi transmitters and multi-antenna WiFi transmitters. The ceiling of the room is organized into a grid pattern, with LEDs positioned at uniform intervals to ensure optimal coverage. The network can serves $K$ users, each equipped with a photodiode and an antenna, randomly distributed within the room. However, for the convenience of analysis, we conducted simulations with $K=3$ users, unless otherwise stated.

\begin{table}[t]
	{
		\begin{center}
			\caption{Simulation Parameters}\label{tbl2}
			\begin{tabular}{ |c|c|c| }
				\hline
				\textbf{Parameter} & \textbf{Value} \\
				\hline
				Room Size & 10m $\times$ 10m $\times$ 5m   \\
				\hline
				Height of Hybrid Transmitters & 4(m)\\
				\hline
				Receiver FOV $(\phi_{FoV})$ & $60^{\circ}$  \\
				\hline
				LED Emission Semiangle $(\psi_{1/2})$ & $70^{\circ}$  \\
				\hline
				Area of Photodiode $(A_D)$ & $1 \: cm^2$  \\
				\hline
				Photodiode Responsivity & 0.54 $A/W$  \\
				\hline
				$\text{LiFi Bandwidth } \left(\text{BW}_{k,s}^\text{LiFi}\right)$  & 20 MHz  \\
				\hline
				$\text{WiFi Bandwidth } \left(\text{BW}_{k,s}^\text{WiFi}\right)$  & 10 MHz  \\
				\hline
				$P_{\text{max}}^{\text{Hybrid}}$ & 38 dbm\\
				\hline
				Refractive Index $(r)$ & 1.5 \\
				\hline
				LiFi Power Amplifier Efficiency $(\eta^\text{LiFi})$ & 0.5 \\
				\hline
				WiFi Power Amplifier Efficiency $ (\eta^{\text{WiFi}}) $ & 0.5 \\
				\hline
				DC Bias $(x^s_{DC,l})$ & $\sqrt{6} \: A$ \\
				\hline
				Maximum Permissible Current $(I_{max,l})$ & 5 $A$\\
				\hline
				\text{Optical Filter Gain } $(G_f)$ & 1 \\
				\hline
				\text{PSD of AWGN in LiFi}  $(\mu_{k,s}^{LiFi})$  & $10^{-19} A^2/Hz$ \\
				\hline
				\text{PSD of AWGN in WiFi } $\left( \mu_{k,s}^\text{WiFi} \right)$ & -174 dbm/Hz\\
				\hline
			\end{tabular}
	\end{center}}
\end{table}	
\begin{figure}[h!]
	%\centering
	\begin{minipage}{0.45\textwidth}
		\centering
		\includegraphics[width=\linewidth]{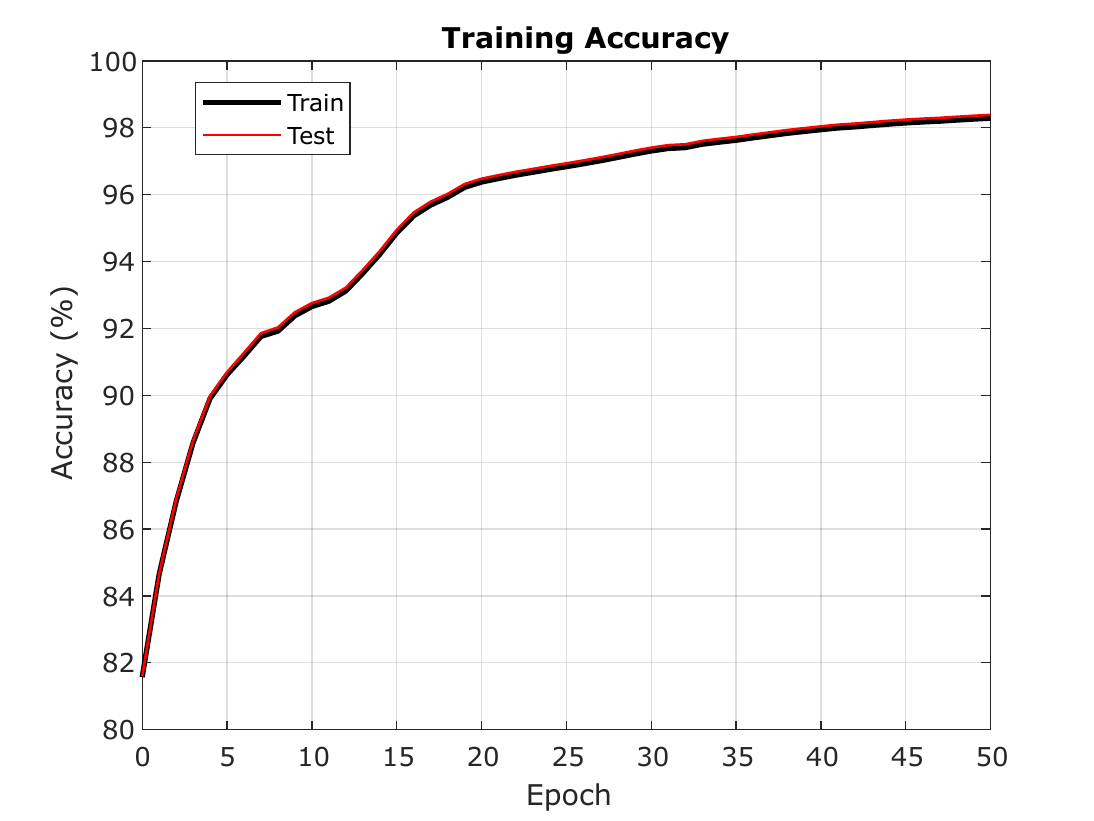}
		\caption{Training Accuracy of the Proposed Model.}
		\label{FIG:1}
	\end{minipage}\hfill
	\begin{minipage}{0.45\textwidth}
		%	\centering
		\includegraphics[width=\linewidth]{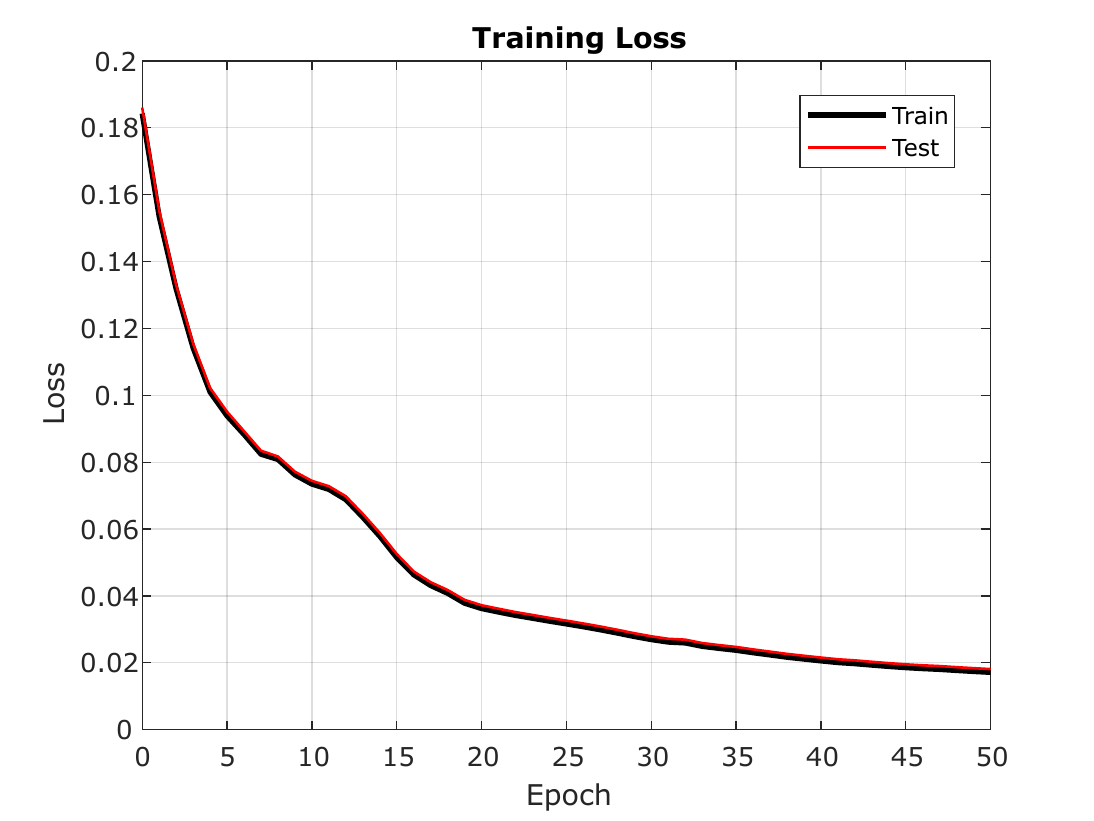}
		\caption{Training Loss of the Proposed Model.}
		\label{FIG:2}
	\end{minipage}
\begin{minipage}{0.45\textwidth}
	\centering
	\includegraphics[width=\linewidth]{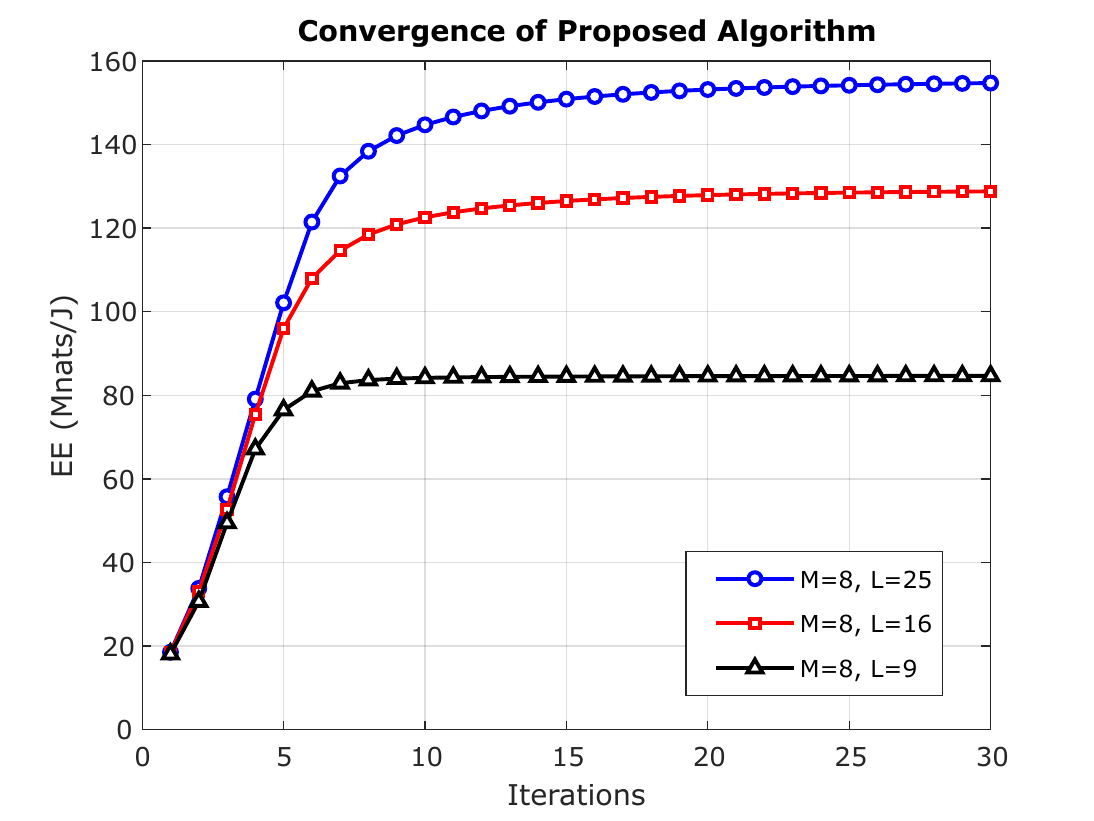}
	\caption{Convergence of the proposed algorithm for maximizing EE as function of precoding vectors}
	\label{FIG:3}
\end{minipage}\hfill
\begin{minipage}{0.45\textwidth}
	%	\centering
	\includegraphics[width=\linewidth]{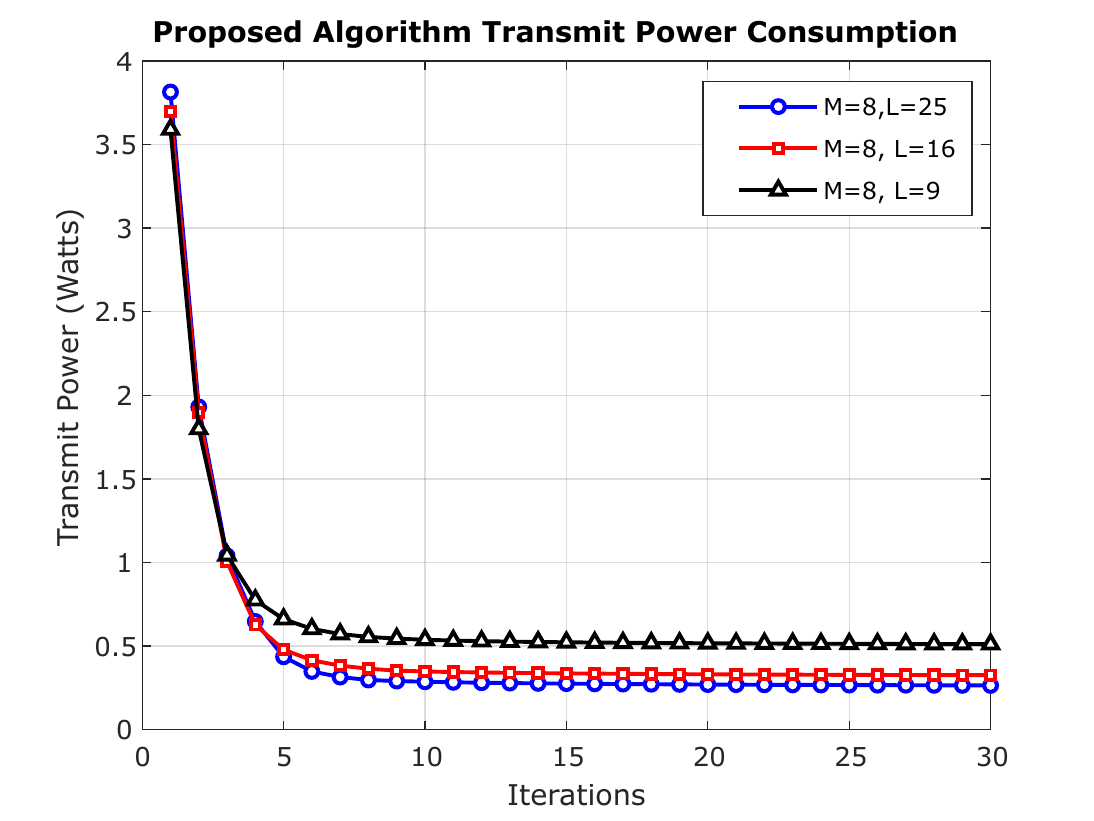}
	\caption{Transmit power minimization based on proposed algorithm.}
	\label{FIG:4}
\end{minipage}
\end{figure}	
	
Our simulation framework employs the mATRIC to dynamically predict the slice type for the user based on monitored KPIs. These KPIs include a range of critical performance indicators, including latency, reliability, supported technologies, use case categories, and temporal aspects of the network infrastructure. To facilitate accurate slice prediction at mATRIC, we train a deep learning model using historical KPI data at universal non-RT RIC. This paper extends the dataset as referenced in \cite{Dataset}, adapting it to hybrid WiFi/LiFi networks. The dataset, includes historical KPI data as features and slice type as label. To prepare this data for modeling, categorical features are converted into a numerical format using one-hot encoding. This process transforms categorical variables into a series of binary columns, each representing a possible category value, thus enabling the model to interpret these features effectively. The label, or output variable, is the slice type, which includes categories such as eMBB, URLLC, and mMTC. The dataset is split into training and testing subsets, with 80\% of the data used for training and 20\% for testing. The deep learning model is implemented and trained using MATLAB's Neural Network Toolbox. Specifically, we use the resilient backpropagation algorithm, with feedforward neural networks configured to have Rectified Linear Unit (ReLU) activation functions in the hidden layers and a softmax function in the output layer for classification. To prevent overfitting, L2 regularization is applied. The performance of the model is assessed by plotting training and validation accuracy and loss over epochs, as shown in Fig. 2 and Fig. 3. These figures provide insights into the model’s capability to accurately predict the appropriate network slice type. The training accuracy of our model is approximately 98 percent for both the training and test datasets, demonstrating its effectiveness in predicting network slices. Additionally, the training loss is relatively low, indicating that the model is performing well in terms of precision and generalization. This trained model can then be utilized by mATRIC for real-time slice selection for users. 

Subsequently, the proposed algorithm optimizes precoding vectors for both LiFi and WiFi technologies to maximize EE with predicted slice types for the users. Our iterative algorithm for maximizing EE shows fast convergence towards the optimal solution, as depicted in Fig. 4. The results demonstrate that our approach effectively addresses the problem, with convergence typically occurring in approximately 10 iterations for 9 LEDs, 15 iterations for 16 LEDs, and 20 iterations for 25 LEDs. This pattern indicates that the number of iterations required for convergence scales with the complexity of the system, specifically the number of LEDs involved. These results were plotted for $M=8$ antennas at the WiFi and varying numbers of LEDs at the LiFi transmitter. The consistent convergence within these iteration ranges highlights the robustness and adaptability of our method across different configurations. Fig. 5 illustrates the power consumption of the Hybrid WiFi/LiFi transmitter during signal transmission for users, as minimized by our proposed algorithm. The figure shows that increasing the number of LEDs at the LiFi transmitter leads to a higher reduction in power consumption compared to using fewer LEDs.

Fig. 6 compares the EE of our proposed precoding algorithm for the hybrid network with ZF and MRT precoding. The comparison is conducted with $M=8$ antennas at the WiFi transmitter and various numbers of LEDs at the LiFi transmitter. The figure clearly shows that our proposed algorithm outperforms the benchmarks in terms of EE. Additionally, it is evident that increasing the number of LEDs at the LiFi transmitter leads to improved EE.

Fig. 7 illustrates the EE of the hybrid network using our proposed algorithm as a function of transmission time $T^t$ for URLLC and mMTC slices, which employ finite block lengths for transmission. As anticipated, the EE improves with increasing transmission time because the system sum-rate increases at a faster rate relative to power consumption. This behavior is evident from the rate expressions for URLLC and mMTC. Additionally, Fig. 7 highlights how EE is influenced by stringent latency requirements. Shorter transmission times, necessitated by strict latency constraints, result in lower EE. For system design purposes, selecting transmission time based on the QoS requirements and desired EE is crucial. These results were plotted for $M=8$ antennas at the WiFi transmitter and varying numbers of LEDs at the LiFi transmitter. Additionally, it can be observed that an increase in the number of LEDs at the LiFi transmitter leads to higher EE with our proposed algorithm.

\begin{figure}[h!]
	%\centering
	\begin{minipage}{0.45\textwidth}
		\centering
		\includegraphics[width=\linewidth]{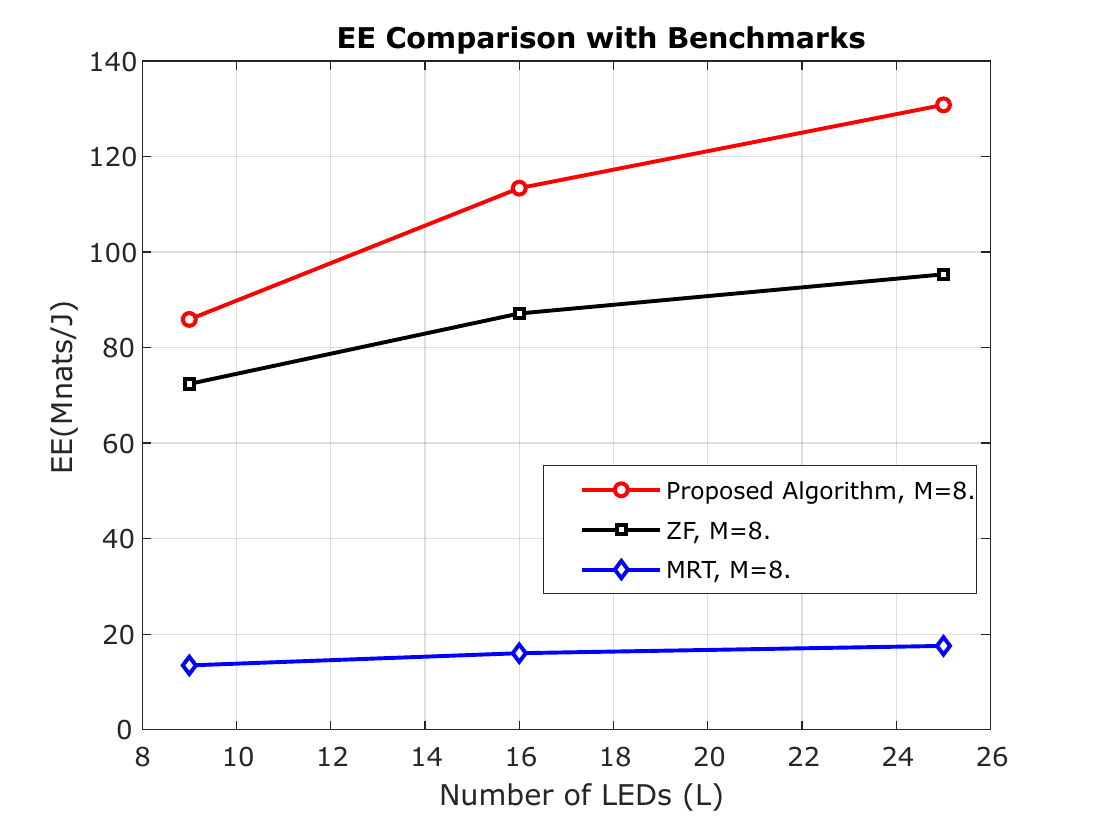}
		\caption{EE comparison of proposed algorithm with benchmarks}
		\label{FIG:5}
	\end{minipage}\hfill
	\begin{minipage}{0.45\textwidth}
		\centering
		\includegraphics[width=\linewidth]{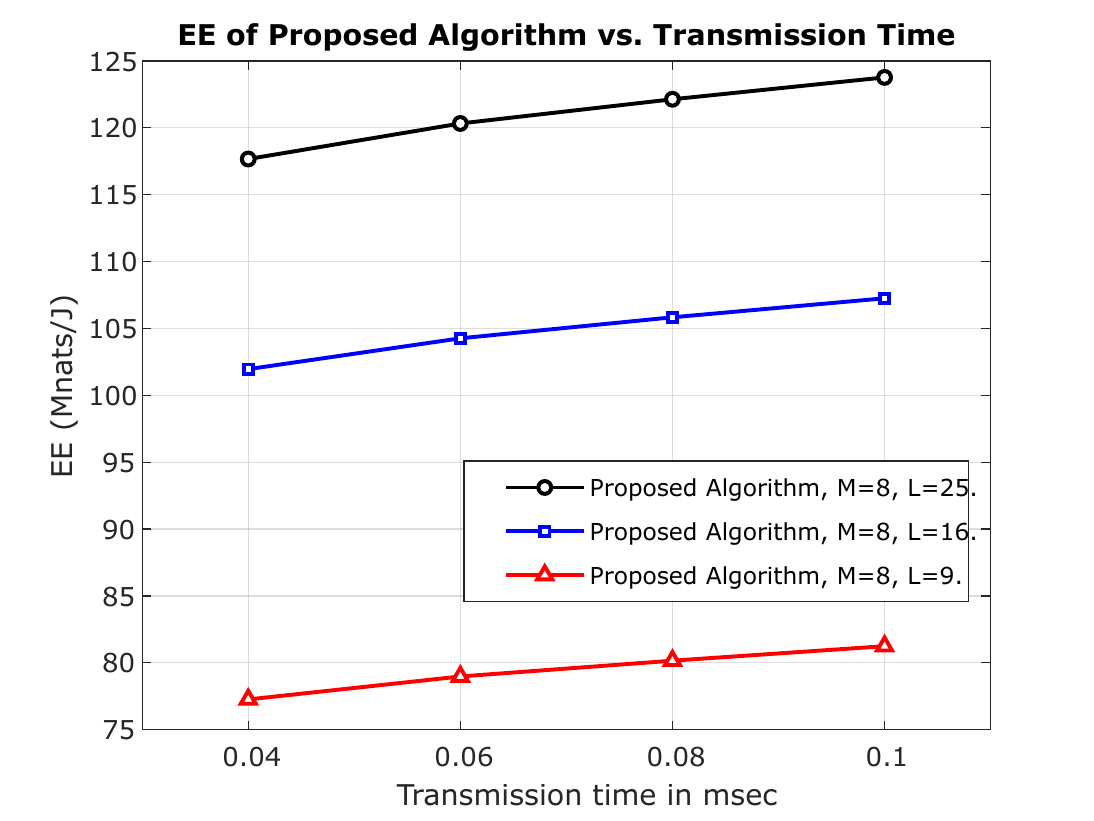}
		\caption{EE of proposed algorithm vs. transmission time}
		\label{FIG:8}
	\end{minipage}
	\begin{minipage}{0.45\textwidth}
		\centering
		\includegraphics[width=\linewidth]{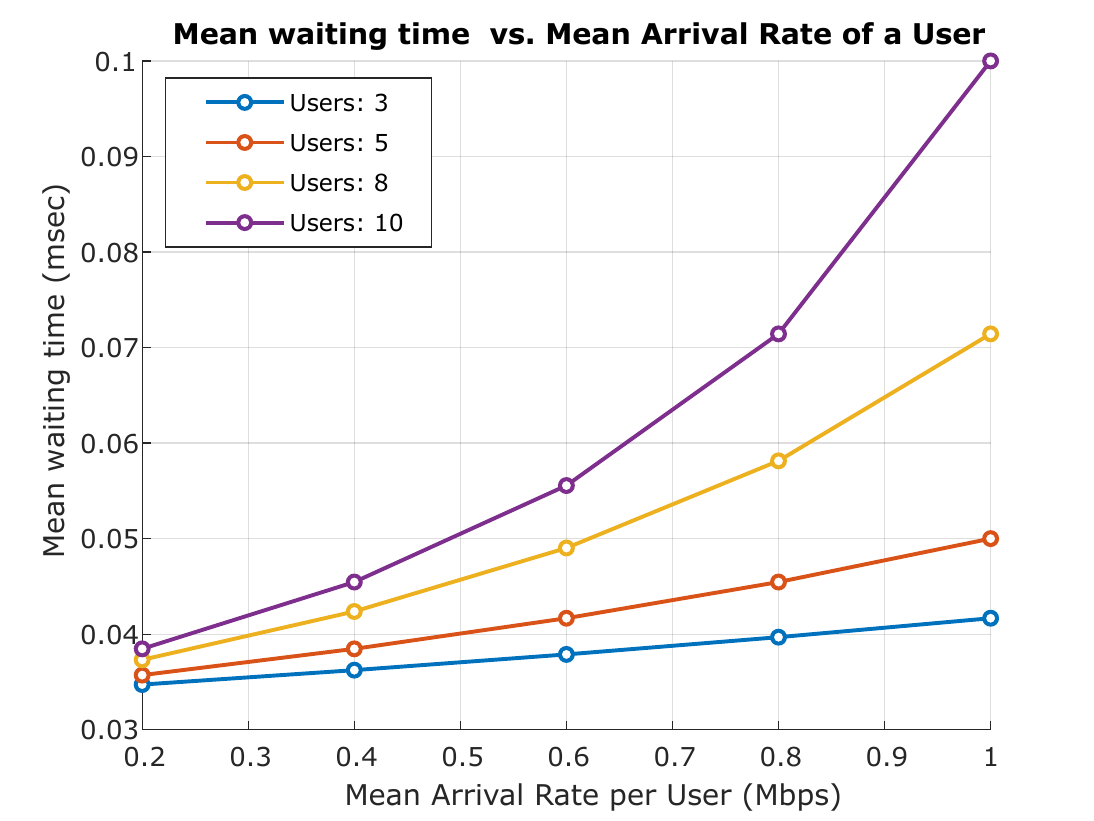}
		\caption{Mean waiting time vs. mean arrival rate of a user}
		\label{FIG:7}
	\end{minipage}\hfill
\begin{minipage}{0.45\textwidth}
	\centering
	\includegraphics[width=\linewidth]{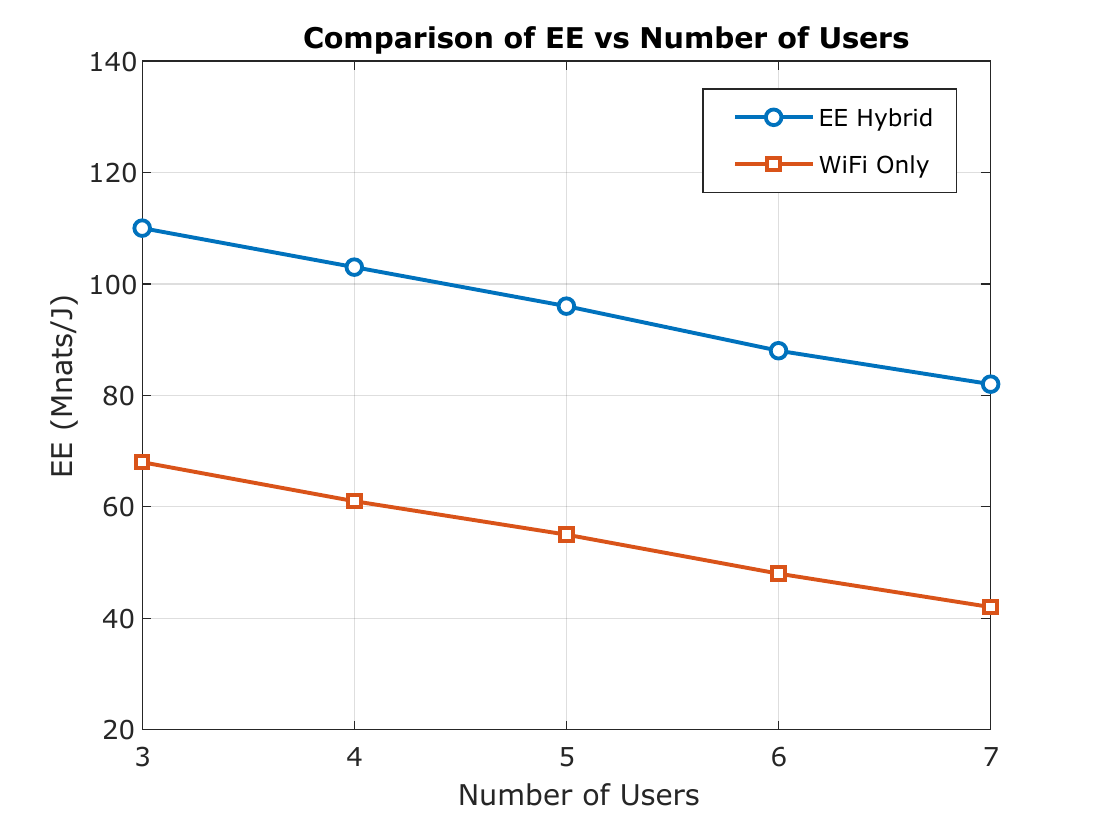}
	\caption{Comparison of EE vs number of users}
	\label{FIG:8}
\end{minipage}
\end{figure}

%The system’s performance is analyzed through a series of tests involving varying numbers of LEDs for LiFi/antennas for WiFi and users, with results presented in the following figures. This setup allows us to assess the effectiveness of our approach in enhancing service differentiation and optimizing energy consumption in a hybrid network environment.
The maximum latency threshold $T^{max}_{k,s}$,  is set based on application needs. For URLLC, it is specified as 1 ms to meet stringent reliability and low latency requirements. For eMBB and mMTC, the thresholds are set at 4 ms and 5 ms, respectively \cite{O-RAN Network slicing}. Within this constraint, transmission time is $T^t_{k,s}= 0.05 \: \text{ms}$, particularly for URLLC, where stringent latency requirements must be met. However, for eMBB and mMTC, longer transmission times can be considered to improve EE as presented in Fig 6. The sum of the channel access delay and backhaul delay is \( T^a_{k,s} + T^b_{k,s} = 0.1 \text{ ms} \), while the combined packet reception and processing delays are \( T^r_{k,s} + T^p_{k,s} = 0.3 \text{ ms} \) for the hybrid WiFi and LiFi network \cite{Latency}. The waiting time \( T^{w}_{s} \) for the users in slice s is modeled as M/M/1 queue model and can be given as \( T^{w}_{s} =\frac{1}{\big( \mu^{Hybrid}_{s}-\lambda_{s}\big)} \), where $\mu^{Hybrid}_{s}$ represents mean service rate of hybrid network for serving the users in slice s, and $\lambda_{s}$ denote the mean arrival data rate for users in slice $s$. It is defined as $
\lambda_s = \sum_{k=1}^{K_s} \alpha_k $, where $\alpha_k$ represents the arrival data rate of user $k$ in the slice. Therefore, the  waiting time of the 
$k^{th}$ user in slice $s$ is \( T^{w}_{k,s} =\frac{1}{\big( \mu^{Hybrid}_{k,s}-\alpha_{k,s}\big)} \). Fig. 8 presents the waiting time experienced by a user in a URLLC slice, plotted against the arrival rate of the users and the total number of users being served. The figure demonstrates that the waiting delay increases with both the mean arrival rate (assuming a constant mean service time) and the number of users in the slice. Moreover, it is observed that if the waiting time for users in the URLLC slice stays below 0.50 ms, the total latency, including all delays, will comfortably meet the 1 ms requirement, thereby ensuring the stringent latency needs of URLLC are fulfilled. Fig. 9 illustrates a comparative analysis of performance between a fully hybrid LiFi/WiFi system and a WiFi-only system, evaluated through iterative optimization of the proposed algorithm. In the hybrid scenario, users benefit from signals transmitted by both WiFi (with \( M = 4 \) antennas) and LiFi (with \( L = 16 \) LEDs). Conversely, the WiFi-only scenario relies solely on WiFi transmitters, equipped with \( M = 4 \) antennas. The results reveal that the hybrid LiFi/WiFi system significantly outperforms the WiFi-only system in terms of EE. This enhancement highlights the advantages of integrating LiFi with WiFi, offering a more effective and resource-efficient solution compared to using WiFi alone. Additionally, the figure demonstrates a decline in EE as the number of users increases in both systems, particularly in scenarios involving mMTC. This reduction is attributed to the rising demand for power resources, which leads to higher power consumption without a corresponding increase in data rates, ultimately diminishing overall EE.

\section{Real-World Applications}
Our proposed algorithm applies to any real-world application where users have different service requirements in the hybrid WiFi/LiFi network. For instance, in health-care settings, various services such as patient monitoring through wearable sensors, emergency messages, remote surgery, and telemedicine including online prescription and consultations have unique demands for data transmission including capacity, latency and reliability, as shown in Fig. 1. Patient monitoring through wearable sensors depends on mMTC service, which can support many wearable sensors. To ensure accurate patient monitoring, the network should send instructions or updates to these sensors with high reliability, moderate capacity and latency. Emergency messages are very critical for fast response scenarios and belong to URLLC service. This service demands extremely low latency and high reliability to offer time-sensitive pieces of information such as messages in life-threatening situations. Remote surgery is also related to URLLC because any interruptions in control signals between the surgeon and robotic instruments can lead to life-threatening consequences. While telemedicine services including online prescription and consultations fall under eMBB service because the transmission of large medical records i.e. high-resolution images and video calls demands higher data rates. Our algorithm is designed to manage these varying demands in the hybrid WiFi/LiFi network simultaneously through network slicing, ensuring that each service receives the necessary resources while maximizing the EE of the entire network.

\section{Conclusion}

	%\section{Conclusion}
	%The conclusion goes here.

	%\section*{Acknowledgments}
	%This should be a simple paragraph before the References to thank those individuals and institutions who have supported your work on this article.
	
	In this paper, we have presented a novel predict-and-optimize algorithm for hybrid WiFi/LiFi networks to address the dual challenges of service differentiation and energy efficiency. By using mATRIC for dynamic network slice prediction based on monitored KPIs and employing a deep learning model trained with the resilient backpropagation algorithm, our framework effectively enables real-time slice selection. This approach not only enhances the adaptability of the network to varying user requirements but also ensures that energy efficiency is maximized. In the optimization phase, we utilized advanced techniques from sequential convex approximation and the inner approximation method to address the complexities associated with non-convex objective functions and constraints. Our novel approximations and iterative algorithm allow for effective convexification, providing a practical means to achieve near-optimal solutions. The performance evaluation through simulations highlights the robustness of our proposed algorithm. The results demonstrate a significant improvement in energy efficiency while maintaining effective service differentiation across hybrid network environments. These findings underscore the potential of our approach to enhance the operational performance of hybrid WiFi/LiFi networks in terms of EE while ensuring service differentiation. This makes our contribution valuable for the ongoing advancements in network optimization and resource management. Future work could focus on further refining the algorithm and exploring its applicability to a broader range of network scenarios and configurations.

	{\appendix [Intial feasible Point for the proposed algorithm]
		
To determine the intial feasible point $\mathbf{f}_{k,s}^\text{WiFi(0)}$, we randomly start from a feasible point for convex power constraint \text{C1-WiFi} and iterate the convex problem
\begin{equation}
	\text{P0:WiFi:} \quad \max_{\mathbf{f}_{k,s}^\text{WiFi}} \Big[\min_{k=1,2,\cdots, K} \overline {S_{k,s}^\text{WiFi}} \Big]
\end{equation}

Subject to:
\begin{equation}
	\begin{aligned}
		& \text{C1-WiFi:} \quad P_{\text{WiFi}} \leq P_{\text{max}}^{\text{WiFi}},
	\end{aligned}
\end{equation}
Similarly for the intial point $\mathbf{f}_{k,s}^\text{LiFi(0)}$, we randomly start from a feasible point for convex power constraint \text{C1-LiFi} and iterate the convex problem
\begin{equation}
	\text{P0:LiFi:} \quad \max_{\mathbf{f}_{k,s}^\text{LiFi}} \Big[\min_{k=1,2,\cdots, K} \overline {S_{k,s}^\text{LiFi}} \Big]
\end{equation}

Subject to:
\begin{equation}
	\begin{aligned}
		& \text{C1-LiFi:} \quad  P_{\text{LiFi}} \leq P_{\text{max}}^{\text{LiFi}},
	\end{aligned}
\end{equation}
The Solution of the above problems can be used as intial feasible point for our main convex problem P3.
}

	%{\appendices
	%\section*{Proof of the First Zonklar Equation}
	%Appendix one text goes here.
	% You can choose not to have a title for an appendix if you want by leaving the argument blank
	%\section*{Proof of the Second Zonklar Equation}
	%Appendix two text goes here.}

	%\newpage
	
	%\section{Biography Section}
	%If you have an EPS/PDF photo (graphicx package needed), extra braces are
	%needed around the contents of the optional argument to biography to prevent
	%the LaTeX parser from getting confused when it sees the complicated
	%$\backslash${\tt{includegraphics}} command within an optional argument. (You can create
	%your own custom macro containing the $\backslash${\tt{includegraphics}} command to make things
	%simpler here.)
	
	%\vspace{11pt}
	
	%\bf{If you include a photo:}\vspace{-33pt}

	%\vspace{11pt}
	
	%\bf{If you will not include a photo:}\vspace{-33pt}
	%\begin{IEEEbiographynophoto}{John Doe}
	%	Use $\backslash${\tt{begin\{IEEEbiographynophoto\}}} and the author name as the argument followed by the biography text.
	%\end{IEEEbiographynophoto}

	\vfill
	
\end{document}